\documentclass[a4paper]{article}
\usepackage{a4wide}
\usepackage{graphicx}      
\usepackage[british]{babel}

\usepackage{tikz}

\usepackage{tikz}
\usetikzlibrary{matrix,positioning}
\usetikzlibrary{shapes,arrows}

\usepackage{cite}
\usepackage{graphicx}
\usepackage{amsmath}
\usepackage{amssymb}

\newcommand{\beq}[1]{\begin{equation}#1\end{equation}}
\newcommand{\bal}[1]{\begin{align}#1\end{align}}
\newcommand{\bpm}[1]{\begin{pmatrix}#1\end{pmatrix}}

\newcommand{\co}{\mbox{co}}

\newcommand{\fci}[3]{\begin{figure}\begin{center} \includegraphics[#1]{{#2}} \caption{{#3}} \end{center}\end{figure}}
\newcommand{\dom}{\ensuremath{\mbox{dom }}}

\usepackage{amsthm}
\newtheorem{thm}{Theorem}

\newtheorem{lem}[thm]{Lemma}
\newenvironment{pf}{\begin{proof}}{\end{proof}}
\theoremstyle{plain}
\newtheorem{dfn}{Definition}
\newtheorem{ass}{Assumption}
\newtheorem{rmk}{Remark}

\begin{document}
\title{Distance function design and Lyapunov techniques for the stability of hybrid trajectories
\thanks{Corresponding author J.~J.~B.~Biemond. {benjamin.biemond@cs.kuleuven.be}. Tel. +32 16 3 27835.
Fax +32 16 3 27996.}
}

\author{J.~J.~Benjamin~Biemond\thanks{Department of Computer Science, KU Leuven, Belgium}, W.~P.~Maurice~H.~Heemels\thanks{Department of Mechanical Engineering, Eindhoven University of Technology, the Netherlands},\\
 Ricardo~G.~Sanfelice\thanks{Department of Computer Engineering,
University of California Santa Cruz, California, U.S.A.} and Nathan van de Wouw \footnotemark[3] \thanks{Department of Civil, Environmental \& Geo- Engineering, University of Minnesota,
Pillsbury Drive SE, U.S.A.}}

\maketitle
%

\begin{abstract}
The comparison between time-varying hybrid trajectories is crucial for tracking, observer design and synchronisation problems for hybrid systems with state-triggered jumps. In this paper, a systematic way of designing an appropriate distance function is proposed that can be used for this purpose.  The so-called ``peaking phenomenon'', which occurs when using the Euclidean distance to compare two hybrid trajectories, is circumvented by taking the hybrid nature of the system explicitly into account in the design of the distance function.
Based on the proposed distance function, we define the stability of a trajectory of a hybrid system with state-triggered jumps and present sufficient Lyapunov-type conditions for stability of a hybrid trajectory. A constructive design method for the distance function is presented for hybrid systems with affine flow and jump maps and a jump set that is a hyperplane. For this case, the mentioned Lyapunov-type stability conditions can be verified using linear matrix conditions. Finally, for this class of systems, we present a tracking controller that asymptotically stabilises a given hybrid reference trajectory, and we illustrate our results with examples.\\[1ex]
\end{abstract}

\section{Introduction\label{secintro}}
Hybrid system models have proven valuable to capture the dynamics of complex systems arising in the domains of mechanical, chemical or electrical engineering, as well as in biological and economical systems, as these models combine continuous-time dynamics with discrete events or jumps \cite{goe_san_12,sch_sch_00,lun_lam_09,hee_sch_10}. While the stability of isolated points or closed sets of hybrid systems is relatively well-understood \cite{goe_san_12,sch_sch_00,lun_lam_09,hee_sch_10}, the stability of time-varying trajectories of these systems received significantly less attention and many issues are presently unsolved. Given the importance of stability of  trajectories in tracking control, observer design and synchronisation problems, it is important to address these open issues.
\\
\emph{Background:}
One of the main complications to study the stability of hybrid trajectories is  the ``peaking phenomenon'' of the Euclidean distance between two trajectories, that can be observed when jump times do not coincide, and the states of two  hybrid trajectories are compared at the same continuous-time instant, as observed in  \cite{lei_wouw_08,men_tor_01} in the framework of measure differential inclusions, in \cite{hee_cam_11} for complementarity systems, and in \cite{san_bie_14,bie_wouw_13} for jump-flow systems.
``Peaking'' of the Euclidean error occurs when two solutions from close initial conditions do not jump at the same time instant. When, before the first jump, the Euclidean error is small, then the Euclidean error approximately equals the jump distance directly after the first jump. A jump of the other solution may again render the Euclidean distance small. As the amplitude of the resulting peak in the Euclidean error cannot be reduced to zero by taking closer initial conditions, trajectories of hybrid systems with state-triggered jumps are generically not asymptotically stable with respect to the Euclidean measure. The latter is even the case when the jump times of both trajectories converge to each other and, consequently, the large Euclidean error occurs only during smaller and smaller time intervals near the jumps, after which periods of flow follow in which both trajectories are close. As this scenario corresponds to desirable behaviour and  therefore should correspond to small errors, it is clear that the Euclidean error is not a good measure in this context.
\\
Only under more severe system assumptions, the standard stability analysis based on the Euclidean error can be employed. Indeed, when jumps of two trajectories are synchronised, then the standard approach can be employed where the difference between both trajectories is required to behave asymptotically stable along trajectories, leading, for example, to successful tracking control approaches presented in  \cite{san_bie_14,lei_wouw_08ijbc}. However, in this paper, we focus our attention to systems with state-triggered jumps, for which, generically, jumps of neighbouring trajectories are not synchronised, such that the ``peaking behaviour'' typically occurs and the comparison of trajectories becomes much more challenging.
\\
The ``peaking'' of the Euclidean error occurs when two states of the system are compared at a given time. Alternatively, the graphs of complete trajectories can be compared. This approach allows to study dynamical characteristics including the continuity of trajectories with respect to initial conditions, as presented in \cite{bro_ara_02,goe_san_12,mor_78}. However, using such a trajectory-based measure, it is hard to formulate constructive conditions (e.g. Lyapunov-based) to guarantee the stability of hybrid trajectories. Therefore, the study of the stability of trajectories in this paper is performed by considering the evolution of a suitably defined distance along the trajectories, therewith necessitating the formulation of a distance function between states of the system.
\\
Focussing on mechanical systems with unilateral position constraints, that include billiard systems, two approaches have been presented in the literature to avoid the ``peaking behaviour''. Firstly, focussing on impacts with non-zero restitution coefficients, it has been observed that during the peaks of the Euclidean error, the two trajectories are far apart, but one trajectory is close to the image of the other that is mirrored in the constraint surface. The Zhuravlev-Ivanov method, cf.\ \cite{bro_99}, describes the trajectory of the real system together with the mirrored images. In this manner, trajectories can cross the constraint by switching from the real to the mirrored state. Hence, the dynamics can be described with a (discontinuous) differential equation without impacts, therewith avoiding the peaking phenomenon. In \cite{for_teel_11_1,for_teel_11_2,for_teel_11tac}, tracking control and observer problems are defined by requiring the asymptotic stability of a set that consists of the real system and the mirrored images. As a second approach, in \cite{gal_men_08, gal_men_12, men_tor_01,mor_bro_10}, the standard Euclidean state error is employed away from the impacts times, while near impacts, only the position error, and no velocity error is considered.
\\
In \cite{bie_wouw_13}, the comparison of trajectories with non-synchronised jumps is facilitated by a distance function that takes the jumping nature of the hybrid system into account, therewith avoiding the ``peaking phenomenon''. In \cite{bie_wouw_13}, we presented sufficient conditions on this distance function, such that stability in this distance function corresponds to an intuitively correct stability notion in the sense that the time mismatch between jumps of trajectories with close initial conditions remains small, and away from the jump times, their states are close. However, no constructive design for this distance function was presented. Only for two examples, such a distance function was proposed in \cite{bie_wouw_13}. Focussing on a class of constrained mechanical systems, a similar distance function was employed in \cite{sch_98} to study continuity of trajectories with respect to initial conditions. Both in \cite{bie_wouw_13} and in \cite{sch_98}, ad-hoc techniques were used to design the distance function.
\\
\emph{Contributions:}
As a first main contribution in the current paper, we present a constructive and general design for the distance function. In order to evaluate this distance function along two different hybrid trajectories, an extended hybrid system is employed, of which each trajectory represents the two original trajectories. This construction results in a \emph{combined} hybrid time domain, and is feasible if both trajectories have a hybrid time domain that is unbounded in the continuous-time direction. We show that when the (global) asymptotic stability is defined with respect to the new distance function, then the proposed  distance function provides an intuitively correct comparison between two hybrid trajectories. As a second main contribution, sufficient conditions for asymptotic stability are presented that rely on Lyapunov functions that may increase during either flow or jump, as long as the Lyapunov function eventually decreases along solutions. For this purpose, maximal or minimal average dwell-time arguments are employed, as proposed in the context of impulsive systems in \cite{hes_lib_08}.
The third main contribution consists of the application of the developed stability theory to tracking control problems for a class of hybrid systems where the jump map is an affine function of the state, the jump set is a hyperplane, and the continuous-time dynamics can be influenced by a bounded control input. This class of systems contains certain models of mechanical systems with unilateral constraints. A piecewise affine tracking control law is designed that achieves asymptotic tracking in the proposed distance measure. This property is proven using a piecewise quadratic Lyapunov function with disconnected sub-level sets, such that the asymptotic stability with respect to the new distance notion can be analysed with computationally tractable matrix relations. Finally, the results of this paper are illustrated with two examples.
\\
\emph{Outline:}
This paper is outlined as follows. First, we present the class of hybrid systems considered in Section~\ref{sechybmodel} and design the distance function in Section~\ref{secdistfunc}. Subsequently, the extended hybrid system is proposed and the stability of trajectories is defined in Section~\ref{secdefstab}. A Lyapunov theorem to study the stability of a hybrid trajectory is presented in Section~\ref{secv}, and a constructive piecewise quadratic Lyapunov function is designed in Section~\ref{secvdes} for a class of hybrid systems with affine jump maps and the jump set contained in a hyperplane. These results are applied to tracking control problems in Section~\ref{secapptrack}. Finally, two illustrative examples are given in Section~\ref{secexamples}, followed by conclusions in Section~\ref{secconcl}.
\\
\emph{Notation:}
Let $\mathbb{N}$ and $\mathbb{N}_{>0}$ denote the set of nonnegative and positive integers, respectively. For a set $X\subset \mathbb{R}^n$, $\partial X$ denotes its boundary and for each $y\in \mathbb{R}^n$, the distance between $y$ and $X$ is $\mbox{dist}(y,X):=\inf_{x\in X}\|x-y\|$. The set $\mathbb{B}\subset \mathbb{R}^n$ is the closed unit ball.
Given a (possibly set-valued) map $F$ with domain of definition $\dom F\subseteq \mathbb{R}^n$ and a set $S\subseteq \dom F$, $F(S)=\{y| y\in F(x),\mbox{ with }x\in S\}$  denotes its image; $F(y)=\emptyset$ for $y\not \in \dom F$,  $F^k(x)$, with $x\in \dom F$, $k\in \mathbb{N}_{> 0}$, denotes $F(F^{k-1}(x))$ and for all $x\in \mathbb{R}^n$, $F^0(x)=\{x\}$. Let $F^{-1}(S)$ denotes its pre-image, namely, $F^{-1}(S)=\{x: F(x)\cap S\neq \emptyset\}$. A set-valued map $F:S\subset \mathbb{R}^n\rightrightarrows \mathbb{R}^n$ is outer semicontinuous if its graph $\{(x,y)\in \mathbb{R}^{n}\times \mathbb{R}^n:\ x\in S,y\in F(x)\}$ is closed, and locally bounded if, for each compact set $\tilde S\subseteq S$, $F(\tilde S)$ is bounded. Using Definition~1.4.11 in \cite{aub_fra_09}, an outer semicontinuous mapping $F: S\rightrightarrows Y$ is \emph{proper} if for every sequence $\{(x_n,y_n)\}_{n\in \mathbb{N}}$ where $y_n\in F(x_n)$ and $y_n$ converges in $Y$, the sequence $\{x_n\}_{n\in \mathbb{N}}$ has a cluster point $\bar x$, i.e., there is a point $(\bar x,\bar y)\in S\times Y$ and a subsequence of $\{x_n,y_n\}_{n\in \mathbb{N}}$ that converges to $(\bar x,\bar y)$. We note that $F$ is proper only if it is outer semicontinuous.
For $n,m\in \mathbb{N}_{>0}$, let $I_n$ and $O_{mn}$ denote the identity matrix and the matrix of zeros of dimension $n\times n$ and $m\times n$, respectively. Given matrices $A,B\in \mathbb{R}^{n\times n}$, $A\prec 0$ and $A\preceq 0$ denote that $A$ is symmetric and negative define or negative semidefinite, respectively. We write $A\preceq B$ and $A\prec B$ when $A-B\preceq 0$ and $A-B\prec 0$, respectively. Similarly, $A\succ B$ denotes $B\prec A$ and $A\succeq B$ denotes $B\preceq A$.
\section{Hybrid system model\label{sechybmodel}}
Consider the hybrid system
\begin{subequations}\label{eqhs}
\bal{
\dot x&\in  F(t,x)&x\in C,\label{eqhsflow}\\
x^+&\in G(x)&x\in D,\label{eqhsj}
}\end{subequations}
with $F:[0,\infty)\times C\rightrightarrows \mathbb{R}^n$ and $G: D\rightrightarrows \mathbb{R}^n$, where $C\subseteq \mathbb{R}^n$ and $D\subseteq \mathbb{R}^n$. We emphasize that the jump map $G$ is independent of the time $t$, which, in the following, will be exploited in the design of the distance function. In contrast to embedding an extra variable with dynamics $\dot t=1$, we prefer to use explicit time-dependency of the flow map $F$, as this allows to study the perturbation of initial conditions without perturbing the initial time. The class of hybrid systems in the form \eqref{eqhs} is quite general and permits modelling systems arising in many relevant applications, including mechanical systems with impacts \cite{goe_san_12} and event-triggered control systems, see e.g. \cite{pos_tab_14}.
\\
To illustrate the ``peaking behaviour'' mentioned in Section~\ref{secintro}, in Fig.~\ref{figpeak}, a reference trajectory $x_d$ and a trajectory $x$ of a hybrid system are shown. The data of this hybrid system is
\bal{
\begin{array}{l}
C=[0,\infty)\times \mathbb{R},\quad D=\{0\}\times (-\infty,0],  \\
F(t,x)=\bpm{x_2&-g+u(t)}^T,\quad G(x)=-x,\label{eqdataexl}
\end{array}
}
with $g=9.81$, such that this systems models the dynamics of a bouncing ball where finite forces can be applied. The reference trajectory $x_d$ is generated by the hybrid system with input $u\equiv 0$, initial condition $x_{d0}=\bpm{0&10}^T$ and initial time $t_0=0$, while the input $u$, that generates the trajectory $x$ from initial condition $x_0=\bpm{0&3}^T$ and initial time $t_0=0$, enforces convergence of $x$ to $x_d$ in the sense that the graphs of both trajectories converge to each other. (In fact, the  tracking control law we propose in Section~\ref{secapptrack} is applied.)
Indeed, the error between the jump times of both trajectories approaches zero over time, and, in addition, away from the jump times the states of both systems approach each other.
However, the Euclidean distance between the trajectories does not converge to zero. This ``peaking phenomenon'' renders the Euclidean distance not appropriate to compare these hybrid trajectories, thereby motivating this study towards systematic techniques to find proper distance functions that do converge to zero in situations as in Fig.~\ref{figpeak}.
\\
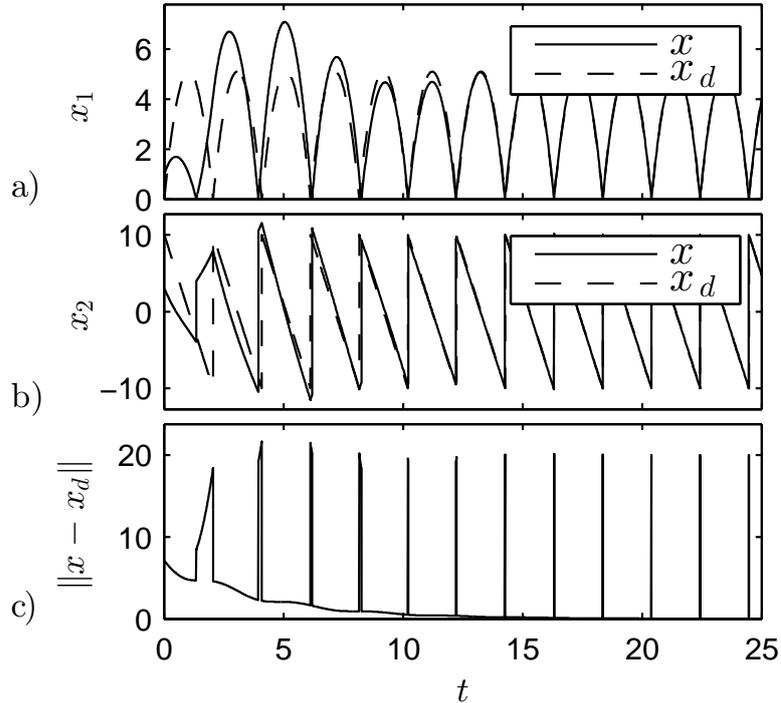
\fci{width=.7\columnwidth}{figxd}{a),b) Projection on the $t$-axis of trajectories $x$ and $x_d$ obtained for the hybrid system with data \eqref{eqdataexl}. c) Euclidean distance function.\label{figpeak}}
We will propose such distance functions for systems \eqref{eqhs} that satisfy the ``hybrid basic conditions'' as defined  for autonomous systems in \cite{goe_san_12}, adapted to allow functions $F(t,x)$ in \eqref{eqhsflow} which depend on $t$. While the conditions in \cite{goe_san_12} are used to ensure robustness and invariance properties, in this paper, the conditions in Assumption~\ref{asssysgen} below are used both to employ Krasovskii-type solutions during flow, and to enable a comparison between trajectories, as will become more clear in Theorem~\ref{thmd} below.
\begin{ass}\label{asssysgen}The data of the hybrid system satisfies
\begin{itemize}
\item $C,D$ are closed subsets of $\mathbb{R}^n$ with $C\cup D\neq \emptyset$;
\item the set-valued mapping $F(t,x)$ is non-empty for all $(t,x)\in [t_0,\infty)\times C$, measurable, and for each bounded closed set $S\subset [t_0,\infty)\times C$, there exists an almost everywhere finite function $m(t)$ such that $\|f\|\leq m(t)$ holds for all $f\in F(t,x)$  and for almost all $(t,x)\in S$;
\item $G:D\rightrightarrows \mathbb{R}^n$ is nonempty, outer semicontinuous and locally bounded.
\end{itemize}
\end{ass}
We consider solutions $\varphi$ to \eqref{eqhs} defined on a hybrid time domain $\dom \varphi\subset [t_0,\infty) \times \mathbb{N}$ as follows, cf.\ \cite{goe_san_12}.
We call a subset $E$ of $[t_0,\infty) \times \mathbb{N}$ a compact hybrid time domain if $E=\bigcup_{j=0}^{J-1}([t_j,t_{j+1}]\times \{j\})$ for some finite sequence $t_0\leq t_1\leq t_2\leq \ldots \leq t_J$. The set $E$ is a hybrid time domain if for all $T,J\in E$, $E\cap ([t_0, T]\times \{0,1,\ldots, J\})$ is a compact hybrid time domain.
Given a hybrid time domain $\dom \varphi$, a hybrid time instant is given as $(t,j)\in \dom \varphi$, where $t$ denotes the ordinary time elapsed and $j$ denotes the number of experienced jumps. The function $\varphi: \dom \varphi\to \mathbb{R}^n$ is a solution of \eqref{eqhs} when jumps satisfy \eqref{eqhsj} and, for fixed $j\in \mathbb{N}$, the function $t\to\varphi(t,j)$ is locally absolutely continuous in $t$ and a solution to \eqref{eqhsflow}. This means $\varphi(t,j)\in D$ and $\varphi(t,j+1)\in G(\varphi(t,j))$ for all $(t,j)\in \dom \varphi $ such that $(t,j+1) \in \dom \varphi$ and
$\varphi(t,j)\in C,$ $\frac{d}{dt}\varphi(t,j)\in \bar F(t,\varphi(t,j))$ for almost all $t\in I_j:=\{t|\ (t,j)\in \dom \varphi\}$ and all $j$ such that $I_j$ has nonempty interior. Herein,
$\bar F(t,x)=\bigcap_{\delta\geq 0}\co \{F(t,(x+\delta \mathbb{B})\cap C)\}$ represents the Krasovskii-type convexification of the vector field which is restricted to $C$, cf.\ \cite{san_goe_08}, where $\co$ denotes the closed convex hull operation. The solution $\varphi$ is said to be complete if $\dom \varphi$ is unbounded.
The hybrid time domain $\dom \varphi$ is called unbounded in $t$-direction when for each $T\geq 0$ there exist a $j$ such that $(T,j)\in\dom\varphi$.
In this paper, we only consider maximal solutions, i.e., solutions $\varphi$ such that there are no solutions $\bar \varphi$ to \eqref{eqhs} with $\varphi(t,j)=\bar \varphi(t,j)$ for all $(t,j)\in \dom \varphi$, and $\dom \bar \varphi$ a hybrid time domain that strictly contains $\dom \varphi$.

\section{Distance function design\label{secdistfunc}}
We will now present a distance function that does not experience the ``peaking behaviour'' that can occur in the Euclidean distance between two trajectories of \eqref{eqhs}, as described in the introduction and illustrated in Fig.~\ref{figpeak}c). We do so for hybrid systems that satisfy the following assumption.
\begin{ass}\label{asssysgtd}The data of the hybrid system \eqref{eqhs} is such that $G$ is a proper function, there is a $k> 0$ for which $G^{k}(D)\cap D= \emptyset$ and every maximal solution of \eqref{eqhs} has a hybrid time domain that is unbounded in $t$-direction.
\end{ass}
We exploit this property in order to define a distance function for the system \eqref{eqhs}. 
\begin{rmk}
Sufficient conditions for the last condition of Assumption~\ref{asssysgtd} can be obtained by an extension of Proposition~2.10 in \cite{goe_san_12} and Lemma 2.7 in \cite{san_goe_07}, which present conditions for completeness and non-Zenoness of trajectories, respectively, towards hybrid systems where the flow dynamics is allowed to be time-dependent, as considered in this paper, see \eqref{eqhsflow}.
\end{rmk}
We now formulate the novel distance function proposed in this paper, where we recall that $G^0(x)$ denotes $\{x\}$ for all $x\in \mathbb{R}^n$.
\begin{dfn}\label{dfnd}
Consider the hybrid system \eqref{eqhs} satisfying Assumption~\ref{asssysgen} and let $\bar k>0$ denote the minimum integer for which Assumption~\ref{asssysgtd} holds. Let the distance function $d:(C\cup D)^2\to \mathbb{R}_{\geq 0}$ be defined by
\beq{\label{eqd}
d(x,y)=\inf_{z\in \mathcal{A}}\left\|\bpm{x\\y}-z\right\|
}
with
\bal{\label{eqa}
\mathcal{A}:=\Big\{\bpm{z_x\\z_y}\in (C\cup D)^2\Big|\ \exists k_1,k_2\in \{0,1,\ldots, \bar k\}, \nonumber\\  G^{k_1}(z_x)\cap G^{k_2}(z_y)\neq \emptyset\Big\}.
}
\end{dfn}
Hence, $d$ vanishes on the set $\mathcal{A}$, which represents all pairs of states $x,y\in C\cup D$ that either are equal or that can jump onto each other by (at most $\bar k$) subsequent jumps characterised by \eqref{eqhsj}.
\\
The following theorem summarises particular properties of the distance function $d$ in Definition~\ref{dfnd}.
\begin{thm}\label{thmdprop}
Consider the hybrid system \eqref{eqhs} satisfying Assumption~\ref{asssysgen} and let $\bar k$ denote the minimum integer for which Assumption~\ref{asssysgtd} holds. The function $d$ in Definition~\ref{dfnd} is continuous and satisfies
\begin{itemize}
 \item [1)] $d(x,y)=0$ if and only if there exist $k_1,k_2\in \{0,1,\ldots,\bar k\}$ such that $G^{k_1}(x)\cap G^{k_2}(y)\neq \emptyset$,
 \item [2)] $\{y\in C\cup D|\ d(x,y)<\beta\}$ is bounded for all $x\in C\cup D$, and all $\beta>0$, and
 \item [3)] $d(x,y)=d(y,x)$, for all $x,y\in C\cup D$.
\end{itemize}\label{thmd}
In addition, the set $\mathcal{A}$ in \eqref{eqa} is closed.
\end{thm}
\begin{pf}
In order to prove 1), we prove that the infimum in \eqref{eqd} is always attained. First, we observe from Assumption~\ref{asssysgen} that $G$ is outer semicontinuous, which directly implies that $G^{-1}$ is outer semicontinuous. In addition, as $G$ is proper according to Assumption~\ref{asssysgtd}, we observe that $G^{-1}$ is locally bounded, cf.\ \cite{aub_fra_09}.
\\
Since the composition $M_1\circ M_2$ of set-valued mappings $M_1$ and $M_2$ is outer semicontinuous and locally bounded when $M_1$ and $M_2$ are outer semicontinuous and locally bounded, we observe that $G^{k_2}$ is outer semicontinuous and locally bounded for all $k_2\in \{0,1,\ldots, \bar k\}$. In addition, reusing this argument, $G^{-k_1}G^{k_2}$ is outer semicontinuous and locally bounded  for all
$k_1,k_2\in \{0,1,\ldots, \bar k\}$.
\\
Note that $\mathcal{A}=\cup_{k_1,k_2\in  \{0,1,\ldots, \bar k\}}A_{k_1k_2}$, with $A_{k_1k_2}:=\{\bpm{x^T&y^T}^T\in (C\cup D)^2|\ y\in G^{-k_1}G^{k_2}(x)\}$, cf.\ \eqref{eqa}. As, for all $k_1,k_2\in \{0,1,\ldots, \bar k\}$, $G^{-k_1}G^{k_2}$ is outer semicontinuous and locally bounded, and $(C\cup D)^2$ is closed, we conclude that each set $A_{k_1k_2}$ is closed. Consequently, we find that the functions $d_{k_1k_2}(x,y):=\mbox{dist}(\bpm{x^T,y^T}^T,A_{k_1k_2})$, for each $k_1,k_2\in \{0,1,\ldots, \bar k\}$, are either continuous functions, or, when $A_{k_1k_2}=\emptyset$, identical to infinity. Since, clearly, $C\cup D\neq \emptyset$ implies $A_{00}$ is nonempty, we observe that $d_{00}(x,y)$ is a continuous and locally bounded function in $C\cup D$. We may write $d(x,y)=\min_{k_1,k_2\in  \{0,1,\ldots, \bar k\}}d_{k_1k_2}(x,y)$, proving that $d$ is continuous. As each set $A_{k_1k_2}$ is closed, $\mathcal{A}$ is closed, such that $d(x,y)=0$ if and only if $\bpm{x^T&y^T}^T\in \mathcal{A}$, proving 1).
\\
We now prove 2) by showing the stronger property that
\beq{\label{prop2infnorm}
Y_\infty(x)\!:\!=\{y\in C\cup D|\ \exists \bpm{z_x\\z_y}\in \mathcal{A}, \|x-z_x\|,\|y-z_y\|\leq \beta\}
} is bounded for each  fixed $x\in C\cup D$ and bounded $\beta>0$. For any $x$, the set $X_\beta^0:=\{w_x|\ \|w_x-x\| \leq  \beta\}$ is compact. Since we have shown above that $G^{-k_1}G^{k_2}$ is outer semicontinuous and locally bounded for all $k_1,k_2\in\{0,1,\ldots, \bar k\}$, we find that the set $G^{-k_1}G^{k_2}(X_\beta^0)$ is compact for all $k_1,k_2\in\{0,1,\ldots, \bar k\}$. As $z_y$ in \eqref{prop2infnorm} has to satisfy $z_y\in G^{-k_2}G^{k_1}(X_\beta^0)$ for some $k_1,k_2=\{0,1,\ldots, \bar k\}$, we have shown that $z_y$ is contained in a bounded set. Hence, we observe that $Y_\infty(x)$ is bounded, which implies 2).
\\
Property 3) directly follows from symmetry of \eqref{eqd} and \eqref{eqa}, which completes the proof.
\end{pf}
\begin{rmk}
Note that the function $d$ in \eqref{eqd} is not a metric, as it does not satisfy the triangle inequality. Namely, if $G$ is set-valued and, for some $x$, $G(x)$ contains two distinct points $y$ and $z$, then $d(x,y)=0$ and $d(x,z)=0$ by Definition~\ref{dfnd}, while $d(y,z)\neq 0$ may still hold in many cases.\\
\end{rmk}
To illustrate that this distance function $d(x,y)$ is non-peaking, in Fig.~\ref{figdnonpeak}, the function $d(x,y)$ is evaluated along the trajectories of Fig.~\ref{figpeak}. While this function is discontinuous in continuous-time $t$ when jumps occur, the function does converge to zero for $t\to \infty$. Hence, the depicted behaviour corresponds to the intuitive observation that the graphs of both trajectories converge towards each other.
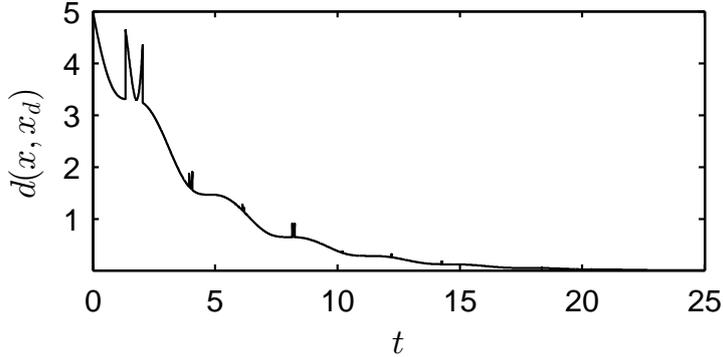
\fci{width=0.7\columnwidth}{figexld}{Distance function $d$ in \eqref{eqd} evaluated along the trajectories shown in Fig.~\ref{figpeak} of the hybrid system with data \eqref{eqdataexl}.\label{figdnonpeak}}
\\
The proposed distance function $d$ in \eqref{eqd} is not contained in the class of functions proposed in \cite{bie_wouw_13}.
Namely, the function $d$ in \eqref{eqd} may not satisfy  $d(x,y)=d(x,g)$ if $y\in D$ and $g\in G(y)$, or $d(x,y)=d(g,y)$ when $x\in D$ and $g\in G(x)$, as was required in \cite{bie_wouw_13}. As another alternative to the distance function in \eqref{eqd}, a more complex distance function design is given in Appendix~\ref{secaltdist}, which in case of $G$ being single-valued and invertible as a function from $D\to G(D)$, ensures that the distance function remains constant during jumps. However, for such distance functions, the set $\{x,y|\ d(x,y)=0\}$ may become undesirably large, in particular when $G$ is not invertible. To allow non-invertible jump maps, we focus on the function $d$ as in \eqref{eqd}.
\begin{rmk}\label{rmkdusef}
To observe that the distance function $d$  in \eqref{eqd} provides an appropriate comparison between two states $x$ and $y$ of \eqref{eqhs}, we observe that a straightforward adaptation of Theorem~1 in \cite{bie_wouw_13} implies that for all $\varepsilon>0$ and all states $x\in C$ with
\beq{
\mbox{dist}(x,D\cup G(D))>\varepsilon\label{distgd},
 }
$d(x,y)<\varepsilon $ implies $\|x-y\|\leq \sqrt 2\varepsilon$ for all $y\in C\cup D$. In many hybrid systems, including models of mechanical systems with impacts, for each solution $\varphi$ the set of times where $x=\varphi(t,j)$ does not satisfy \eqref{distgd} becomes very small when $\varepsilon$ is reduced. Hence, if in these systems $d(x,y)$ is sufficiently small, then $\|x-y\|$ will be small away from some ``peaks'', which can only occur in small time intervals.
\end{rmk}

\section{Stability of hybrid trajectories\label{secdefstab}}
We now evaluate the distance function $d$ along trajectories $\varphi_x(t,j)$, $\varphi_y(t,j)$ of \eqref{eqhs}. In order to enable the comparison of the states of two trajectories in terms of the distance $d$, inspired by the approaches in \cite{san_bie_14} for time-triggered jumps and in \cite{bie_wouw_13} for state-triggered jumps, we introduce an extended hybrid system with state $q\in (C\cup D)^2$, such that a combined hybrid time domain is created. The first and second collection of $n$ components of $\varphi_q(t,j)$, with $\varphi_q$ being the solution to the extended hybrid system, contain a representation of the trajectories $\varphi_x(t,j)$ and $\varphi_y(t,j)$ of \eqref{eqhs} on a `combined' hybrid  time domain.
\\
For this purpose, we construct an extended hybrid system with state $q=\bpm{x^T&y^T}^T\in (C\cup D)^2$, continuous dynamics
\begin{subequations}
\label{eqhse}
\bal{
\dot q=&\bpm{\dot x\\\dot y}\in F_e(t,q):=\bpm{F(t,x)\\F(t,y)} &  \bpm{x\\y}\in C_{e}:=C^2,
}
and jumps characterised by
\bal{
&q^+\!=\!\bpm{x^+\\y^+}\!=\!G_e(q)\!:=\!\begin{cases}
\bpm{G(x)\\ y}\mbox{ if }  x\in D, y\in C\setminus D&\\
\bpm{x\\ G(y)}\mbox{ if } x\in C\setminus D, y\in D&\\
\left\{\bpm{G(x)\\ y},\bpm{x\\ G(y)}\right\}\\\hspace{2.7cm}\mbox{ if }x,y\in D&
\end{cases}\nonumber
\\&\mbox{for }q\in D_{e}\!:\!=\left\{q\!=\!\bpm{x\\y}\in (C\cup D)^2\bigg|\ x\in D\lor y\in D\right\}.
}
\end{subequations}
Given the initial conditions $\varphi_x(t_0,0)$ and $\varphi_y(t_0,0)$ at initial time $(t_0,0)$ for the individual trajectories $\varphi_x,\varphi_y$, respectively, we select the initial condition $\varphi_q(t_0,0)=\bpm{\varphi_x^T(t_0,0)&\varphi_y^T(t_0,0)}^T.$
\\
Solutions of this extended system generate a combined hybrid time domain. This allows to compare two trajectories of the hybrid system at every hybrid time instant $(t,j)\in \dom \varphi_q$. Hereto, let
\beq{\label{xybars}
\begin{matrix}
\bar \varphi_x(t,j)&:=\bpm{I_n&O_{nn}}\varphi_q(t,j),\\
\bar \varphi_y(t,j)&:=\bpm{O_{nn}&I_n}\varphi_q(t,j),
\end{matrix}
}
such that at every time instant $(t,j)\in \dom \varphi_q$, the distance $d(\bar \varphi_x(t,j),\bar \varphi_y(t,j))$ can indeed be evaluated.

\begin{rmk}\label{rmkrepar}
We note that when one of the two trajectories $\varphi_x(t,j)$ and $\varphi_y(t,j)$ has a time domain that is bounded in $t$-direction, then this extended system does not represent both trajectories completely, cf.\ Assumption~\ref{asssysgtd}. Namely, if a trajectory of \eqref{eqhs} (say, the trajectory $\varphi_x$), has a time domain that is bounded in $t$-direction, such that $t\leq T,$ for all $(t,j)\in \dom \varphi_x$, then $s\leq T$ for all $(s,j)\in \dom \varphi_q$ as well, with $\varphi_q(t,j)$ the corresponding solution to \eqref{eqhse}. To see this, note that when $\dom \varphi_x$ is bounded in $t$-direction, then $\varphi_x$ leaves $C\cup D$, has a finite escape time, or has an accumulation of jumps (i.e., experiences Zeno-behaviour), cf.\ \cite[Proposition~2.10]{goe_san_12}. By the construction of \eqref{eqhse}, we observe that $\varphi_q$ also leaves $C_{e}\cup D_{e}$, has a finite escape time, or has an accumulation of jumps. If $\dom \varphi_y$ contains a hybrid time $(t,j)$ with $t>T$, then the trajectory $\varphi_y$ at this time instant is not captured in the dynamics of \eqref{eqhse}.
\\
If both trajectories are unbounded in $t$-direction, then the functions $\bar \varphi_x,\bar \varphi_y$ in \eqref{xybars} are reparameterisations of trajectories $\varphi_x,\varphi_y$ of \eqref{eqhs}. To be precise, there exist non-decreasing functions $j_x,j_y: \mathbb{N}\to \mathbb{N}$ such that $\bar \varphi_x(t,j)=\varphi_x(t,j_x(j))$ and $\bar \varphi_y(t,j)=\varphi_y(t,j_y(j))$, for all $(t,j)\in \dom \varphi_q$.
\end{rmk}
We now employ this combined hybrid time domain and the distance function \eqref{eqd} in order to define the stability of trajectories for hybrid systems.
\\
The distance function defined in \eqref{eqd} allows to compare different points $x,y\in C\cup D$ while taking the jumping nature of the hybrid system \eqref{eqhs} into account. We will now define the stability of trajectories for hybrid systems analoguous to the  definition for ordinary differential equations, cf.\ \cite{leo_08}, by replacing the standard (often Euclidean) metric by the distance function $d$ in \eqref{eqd}. In this manner, we obtain a stability notion that allows ``peaking'' of the standard metric, cf.\ \cite{bie_wouw_13}.\\
Given a trajectory $\varphi_x$ of \eqref{eqhs}, we say that a trajectory $\bpm{\bar \varphi_x^T& \bar \varphi_y^T}^T$ of \eqref{eqhse} represents $\varphi_x$ in the first $n$ states when $\bar \varphi_x$ is a reparameterisation of $\varphi_x$ as in Remark~\ref{rmkrepar}. Clearly, any trajectory to \eqref{eqhse} represents $\varphi_x$ in the first $n$ states when both $\bar \varphi_x(t_0,0)=\varphi(t_0,0)$ holds and this initial condition has a unique solution to \eqref{eqhs}, as considered in \cite{bie_wouw_13}.
%
\begin{dfn}\label{dfnstabtraj}
Consider a hybrid system \eqref{eqhs} satisfying Assumption~\ref{asssysgtd} and let $d$ be given in \eqref{eqd}. The trajectory $\varphi_x$ of \eqref{eqhs} is called \emph{stable with respect to $d$} if for all $\epsilon>0$ there exists a $\delta(\epsilon)>0$ such that for every initial condition $\varphi_y(t_0,0)$ satisfying $d(\varphi_x(t_0,0),\varphi_y(t_0,0))\leq \delta(\epsilon)$, it holds that
\beq{
d(\bar \varphi_x(t,j),\bar \varphi_y(t,j))< \epsilon\mbox{ for all }(t,j)\in \dom \varphi_q,
}
with $\varphi_q(t,j)=\bpm{\bar \varphi_x(t,j)\\\bar \varphi_y(t,j)}$ being any trajectory of the combined system \eqref{eqhse} with initial condition\\
$\bpm{\varphi_x(t_0,0)^T&\varphi_y(t_0,0)^T}^T$ that represents $\varphi_x$ in the first $n$ states, and is called \emph{asymptotically stable with respect to $d$}  if $\delta$ can be selected such that, in addition,
\beq{
\lim_{t+j\to \infty} d(\bar \varphi_x(t,j),\bar \varphi_y(t,j))=0.\label{eqdconv}
}
When the trajectory $\varphi_x$ is asymptotically stable with respect to $d$ and  \eqref{eqdconv} holds for all solutions $\varphi_q$ to \eqref{eqhse} representing $\varphi_x$ in the first $n$ states, then the trajectory $\varphi_x$ is called
globally asymptotically stable with respect to $d$.
\end{dfn}
\begin{rmk}
This notion of the stability of a trajectory $\varphi_x$ corresponds to an intuitive stability notion if convergence of $d(\bar \varphi_x(t,j),\bar \varphi_y(t,j))$ to zero implies that, firstly, the Euclidean distance $\|\bar \varphi_x(t,j)-\bar \varphi_y(t,j)\|$ converges to zero apart from some ``peaks'' near the jumps and, secondly, that the time mismatch of the jumps, which coincides with the time duration of these peaks, converges to zero.
\\
Indeed, we can identify an important class of trajectories $\varphi_x$ for which these implications holds. Given a trajectory $\varphi_x$, let $T_{nd}^\tau(\epsilon):=\{t\geq \tau|\  \mbox{dist}(\varphi_x(t,\bar j
),D\cup G(D))\leq \epsilon,\mbox{ for all } \bar j \mbox{ with }(t,\bar j)\in \dom \varphi_x\}$ for $\epsilon>0$. Consider the class of trajectories $\varphi_x$ for which the Lebesgue measure of $T_{nd}^\tau(\epsilon)$, for sufficiently large $\tau$, tends to zero for $\epsilon\to 0$ (for example, this class contains periodic trajectories that visit $D\cup G(D)$ finitely many times each period, or when the trajectory remains bounded away from $D\cap G(D)$). When such a trajectory is asymptotically stable with respect to $d$, we observe that for any trajectory $\varphi_y$ with $d(\varphi_x(t_0,0),\varphi_y(t_0,0))$ sufficiently small, after some transient, $d(\bar\varphi_x(t,j),\bar\varphi_y(t,j))<\epsilon$, with  small $\epsilon>0$. According to Remark~\ref{rmkdusef}, for every $\epsilon>0$ there exists $T>0$ such that $\|\bar \varphi_x(t,j)-\bar \varphi_y(t,j)\|<\sqrt 2\epsilon,$ for all $(t,j)\in \dom \varphi_q, t\not \in T_{nd}^T(\epsilon)$.
Hence, since the measure of $T_{nd}^\tau(\epsilon)$ tends to zero for $\varepsilon\to 0$, the duration of possible ``peaking'' in the Euclidean error tends to zero, and away from these peaks, $\|\bar \varphi_x(t,j)-\bar \varphi_y(t,j)\|$ tends to zero.
\end{rmk}
%
%
%
%

\section{Lyapunov conditions for stability of trajectories with respect to $d$\label{secv}}
Now, we present sufficient conditions for stability of a trajectory of the system \eqref{eqhs} in the sense of Definition~\ref{dfnstabtraj}, that are based on the existence of an appropriate Lyapunov function. In order to allow the Lyapunov function to increase during flow, and decrease during jumps, or vice versa, the following definitions of minimal and maximal average inter-jump time are adapted from \cite{san_bie_14}. 
\begin{dfn}$ $\label{defdwelltime}
\begin{itemize}
\item A hybrid time domain $E$ is said to have minimal average inter-jump time $\tau > 0$ if there exists $N_0 > 0$ such that for all $(t, j) \in E$ and all $(T,J)\in E$ where $T+J\geq t+j$,
\end{itemize}
\beq{
J-j \leq  N_0 +\frac{T-t}{\tau}.\label{defsavg}}
\begin{itemize}
\item
A hybrid time domain $E$ is said to have maximal average inter-jump time
$\tau > 0$, if there exists $N_0 > 0$ such that for all $(t, j) \in E$ and all $(T,J)\in E$ where $T+J\geq t+j$,

\end{itemize}
\beq{J-j \geq \frac{T-t}{\tau}-N_0.\label{defsravg}
}
We say that a hybrid trajectory $\varphi_q$ has a minimal or maximal average inter-jump time if $\dom \varphi_q$ has a minimal or maximal average inter-jump time, respectively.
\end{dfn}
The following theorem presents Lyapunov-based sufficient conditions for the stability of a trajectory of \eqref{eqhs}. When the trajectories of \eqref{eqhs} have a minimal or maximal average inter-jump time, the requirements on the data of \eqref{eqhs} is less restrictive than in the generic case. As we are interested in the stability of a trajectory, these conditions are imposed only near this trajectory. For this purpose, we recall that given a function function $V:\mathbb{R}^{2n}\to \mathbb{R}_{\geq 0}$ and scalar $v_L>0$, $V^{-1}([0,v_L])$ denotes $\{q\in \mathbb{R}^{2n}|\ V(q)\in [0,v_L]\}$
\begin{thm}\label{thmgen}
Consider a hybrid system \eqref{eqhs} satisfying Assumptions~\ref{asssysgen} and \ref{asssysgtd}. Let $d$ be given in \eqref{eqd}. The trajectory $\varphi_x$ of system \eqref{eqhs} is asymptotically stable with respect to $d$ if there exist a continuous function $V:\mathbb{R}^n\times \mathbb{R}^n\to \mathbb{R}_{\geq 0}$, $\mathcal{K}_\infty$-functions $\alpha_1,\alpha_2$, a scalar $v_L>0$ and scalars $\lambda_c,\lambda_d$ such that $V$ is continuously differentiable on an open domain containing $V_L:=V^{-1}([0,v_L])$ and, for all $(t,j)\in \dom \varphi_x$, it holds that
\bal{
&\alpha_1(d(\varphi_x(t,j),y))\leq V(\varphi_x(t,j),y) \leq \alpha_2(d(\varphi_x(t,j),y)),\nonumber \\&\hspace{4mm}\mbox{for all } y \mbox{ such that} \bpm{\varphi_x(t,j)\\y}\in C_{e}\cup D_{e},\label{eqvbnds}\\
&V(g)\leq e^{\lambda_d} V(q), \mbox{ for all } g\in G_e(q),\nonumber \\&\hspace{4mm}\mbox{and all }y\mbox{ such that }q=\bpm{\varphi_x(t,j)\\y}\in D_{e}\cap V_L\label{eqvjump},\\
&\Big\langle\! \left.\frac{\partial V}{\partial q}\right|_{q},f\!\Big \rangle \leq \lambda_c V(\varphi_x(t,j),y)\mbox{ for all } f\in F_e(t,q)\nonumber \\&\hspace{4mm}\mbox{and all } y\mbox{ such that }q=\bpm{\varphi_x(t,j)\\y}\in C_{e}\cap V_L\label{eqvflow},
}
and at least one of the following conditions are satisfied:
 \begin{enumerate}
 \item[1)] $\lambda_c<0,\lambda_d\leq 0$;
  \item[2)] all trajectories of \eqref{eqhs} have minimal average inter-jump time
$2\tau > 0$,$\lambda_c\leq 0$ and $\lambda_d+\lambda_c \tau<0$;
  \item[3)] all trajectories of \eqref{eqhs} have maximal average inter-jump time
$2\tau > 0$, $\lambda_d\leq 0$ and $\lambda_d+ \lambda_c \tau<0$.
\end{enumerate}
When \eqref{eqvjump} and \eqref{eqvflow} hold for all $y$ such that $q=\bpm{\varphi_x(t,j)^T&y^T}^T\in D_{e}$ and $C_{e}$, respectively, then $\varphi_x$ is globally asymptotically stable with respect to $d$.
\end{thm}
\begin{pf}
We restrict our attention to trajectories $\varphi_q$ to \eqref{eqhse} that represent $\varphi_x$ in the first $n$ states. These trajectories always exist, which follows from the comparison of \eqref{eqhs} and \eqref{eqhse} and the observation that $\varphi_x$ is a trajectory to \eqref{eqhs}. The observation that $\bar \varphi_y$ given in \eqref{xybars} is a reparameterisation of a trajectory $\varphi_y$ for \eqref{eqhs}, and both $\varphi_x$ and $\varphi_y$ are unbounded in $t$-direction by Assumption~\ref{asssysgtd}, proves that the trajectory $\varphi_q$ is unbounded in $t$-direction.
\\
We first prove that $V(\varphi_q(t,j))< v_L$ for all $(t,j)\in \dom \varphi_q$ and all trajectories $\varphi_q$ of \eqref{eqhse} if $\bar kV(\varphi_q(t_0,0))<v_L$, where $\bar k$ is chosen as $\bar k=1$ if 1) holds, $\bar k=e^{\lambda_d N_0}$ if 2) holds and $\lambda_d\geq 0$, and $\bar k=e^{\lambda_c N_0\tau}$ if 3) holds and $\lambda_c\geq 0$, with $N_0$ given in Definition~\ref{defdwelltime}. Observe that if all trajectories of \eqref{eqhs} have a minimal or maximal average inter-jump time $2\tau$, then \eqref{eqhse} has minimal or maximal average interjump time $\tau$.
\\
To prove that the values of $\bar k$ defined above are appropriate, for the sake of contradiction, suppose that $\bar kV(\varphi_q(t_0,0))<v_L$ and there exists a time $(t_0+\bar T,\bar J)\in \dom \varphi_q$, $\bar T,\bar J\geq 0$, such that $V(\varphi_q(t_0+\bar T,\bar J))\geq v_L$. Hence, there exist $T\leq \bar T$ and $J\leq \bar J$ such that $(t_0+T,J)\in \dom \varphi_q$, and
\beq{
V(\varphi_q(t_0+T,J))\geq v_L,\label{contrTJ}
}
but $V(\varphi_q(t,j))< v_L$ for all $(t,j)\in R:=\{(t,j)\in \dom \varphi_q|\ t<t_0+T\lor j<J\}$.
\\
Given the inequalities \eqref{eqvbnds}-\eqref{eqvflow} and the fact that $\varphi_q$ represents $\varphi_x$ in the first $n$ states, we observe that we may replace $\varphi_x$ with $\bar \varphi_x$ and $(t,j)\in \dom \varphi_x$ with $(t,j)\in \dom \varphi_q$ in these inequalities, as $\bar \varphi_x(t,j)=\varphi_x(t,j_x(j)),$ for all $(t,j)\in \dom \varphi_q$, cf.\ Remark~\ref{rmkrepar}. Hence, $V(g)\leq e^{\lambda_d}V(\varphi_q(t,j))$ and $\Big\langle\! \left.\frac{\partial V}{\partial q}\right|_{\varphi_q(t,j)},f\!\Big \rangle \leq \lambda_c V(\varphi_q(t,j))$ hold for all $(t,j)\in R$, $f\in F_e(t,\varphi_q(t,j))$ and $g\in G_e(\varphi_q(t,j))$.
\\
Analogue to \cite{san_bie_14}, we
study the function $t,j\mapsto w(t,j):=V(\bar \varphi_x(t,j),\bar \varphi_y(t,j))$ along the given solution $\varphi_q$ over the time domain $R$ and we introduce scalars $\{t_j\}$ such that $R=\bigcup_j ([t_j,t_{j+1}]\times \{j\})$. As, for each $j$, $\bar \varphi_x,\bar \varphi_y$ are absolutely continuous in $t$ in the time interval $[t_j,t_{j+1}]\times \{j\}$, $w(t,j)$ is absolutely continuous in $t$ as well. Evaluating $\dot w(t,j)=\tfrac{\partial V}{\partial q} \bpm{\frac{d \bar \varphi_x^T(t,j)}{dt}&\frac{d \bar \varphi_y^T(t,j)}{dt}}^T=\tfrac{\partial V}{\partial q} f$ for some $f\in F_e(t,\bpm{\bar \varphi_x^T(t,j)&\bar \varphi_y^T(t,j)}^T)$, we find with  \eqref{eqvflow} that $\dot w(t,j)\leq \lambda_c w(t,j)$. With the comparison lemma, \cite[Lemma 3.4]{kha_02}, we find $w(t_{j+1},j)=e^{\lambda_c (t_{j+1}-t_j)} w(t_{j},j)$ for all $j$. For a subsequent jump, \eqref{eqvjump} yields $w(t_{j+1},j+1)=e^{\lambda_d}w(t_{j+1},j)$. Applying this result repetitively, we find
\beq{\label{eqvtj}
w(t_0+T,j)=V(\varphi_q(t_0+T,J))\leq e^{\lambda_c T+\lambda_d J}V(\varphi_q(t_0,0)).
}
If case 1) of the theorem holds, we directly observe $V(\varphi_q(t_0+T,J))\leq V(\varphi_q(t_0,0))$, contradicting \eqref{contrTJ}. If $\lambda_d\geq 0$ and case 2) holds, then the definition of minimal average inter-jump time yields $\lambda_c T+\lambda_d J\leq \frac{T}\tau (\lambda_c \tau+\lambda_d)+\lambda_d N_0\leq \lambda_d N_0$, such that with \eqref{eqvtj} we find $V(\varphi_q(t_0+T,J))\leq \bar k V(\varphi_q(t_0,0))<v_L$, contradicting \eqref{contrTJ}.
If $\lambda_c\geq 0$ and case 3) holds, then applying the definition of maximal average inter-jump time, we observe that $\lambda_c T+\lambda_d J\leq (\lambda_d+\lambda_c \tau) J+\tau N_0 \lambda_c\leq \lambda_c\tau N_0$. Substituting this inequality in \eqref{eqvtj} we find $V(\varphi_q(t_0+T,J))\leq \bar k V(\varphi_q(t_0,0))<v_L$, contradicting \eqref{contrTJ}. A contradiction has been obtained in all three cases, proving that $\bar k V(\varphi_q(t_0,0))<v_L$  implies $\varphi_q(t,j)\in V_L$ for all $(t,j)\in \dom \varphi_q$.
Hence, $V(\varphi_q(t_0,0))\leq \frac {v_L} {\bar k}$ implies that, for all $(t,j)\in \dom \varphi_q$,  $V(\varphi_q(t_0+t,j))\leq e^{\lambda_c t+\lambda_d j}V(\varphi_q(t_0,0))$.
\\
Assumption~\ref{asssysgtd} states that all trajectories of \eqref{eqhs} are unbounded in $t$-direction, which implies $G(D)\subseteq C\cup D$. Hence, we find $\varphi_q(t_0+t,j)\in C_e\cup D_e$ and we can use \eqref{eqvbnds}. Consequently, we find $d(\varphi_q(t_0+t,j))\leq \alpha_1^{-1}(e^{\lambda_c t+\lambda_d j}\alpha_2(d(\varphi_q(t_0,0))))$. With the inequalities for $\lambda_c t+\lambda_d j$ derived above, we conclude that in any of the three cases of the theorem,
 $d(\varphi_q(t_0+t,j))\leq \alpha_1^{-1}(\bar k\alpha_2(d(\varphi_q(t_0,0))))$, proving stability with respect to $d$. Again using the mentioned inequalities, we observe that $\lambda_c t+\lambda_d j\to -\infty$ along the solutions (this limit can be used since all trajectories are unbounded in $t$-direction, cf.\ Assumption~\ref{asssysgtd}), such that $d(\varphi_q(t_0+t,j))\to 0$. This proves asymptotic stability.
\\
When \eqref{eqvjump} and \eqref{eqvflow} hold for all $y$ such that\\ $\bpm{\varphi_x(t,j)^T&y^T}^T\in C_{e}\cup D_{e}$, then the  upper bounds on $d(\varphi_q(t_0+t,j))$ prove global asymptotic stability.
 \end{pf}

We note that the conditions as presented in this lemma can be extended to time-varying Lyapunov functions $V$. Furthermore, as all trajectories of \eqref{eqhse} are complete by Assumption~\ref{asssysgtd}, it directly follows that $G(D_e)\subset D_e\cup C_e$, such that \eqref{eqvbnds} implies that $V$ is defined on $C_e\cup D_e\cup G(D_e)$.
\begin{rmk}
The Lyapunov function $V$ in this theorem is closely related to the Lyapunov functions used for \emph{incremental stability}, see e.g. \cite{ang_02,ruf_wouw_13} for ordinary differential equations and \cite{li_phi_14} for hybrid systems where incremental stability is defined with respect to the Euclidean distance, and Lyapunov functions in \cite{zam_wouw_13} where incremental stability with respect to non-Euclidean distance functions is investigated for ordinary differential equations. In fact, if the conditions of Theorem~\ref{thmgen} hold for any solution $\varphi_x(t,j)$ of \eqref{eqhs}, then they imply an incremental stability property with respect to the distance $d$. Sufficient conditions for this more restrictive system property are attained by replacing $\varphi_x(t,j)$ in \eqref{eqvbnds}-\eqref{eqvflow} with $x$ and requiring the conditions to hold for all $\bpm{x^T&y^T}^T\in C_e\cup D_e$.
\end{rmk}
For the specific class of hybrid systems with a jump map that is an affine function of the state, and a jump set that is a subset of a hyperplane, in the following section, a piecewise quadratic Lyapunov function is presented which, locally, satisfies the requirements \eqref{eqvbnds} and \eqref{eqvjump} by design. Hereby, we provide a constructive Lyapunov-based approach for (local) stability analysis of trajectories for this class of hybrid systems.

\section{Constructive Lyapunov function design for hybrid systems with affine jump map\label{secvdes}}
We now focus on the class of hybrid systems that have single-valued, affine and invertible jump maps and have jump sets characterised by a hyperplane.
In addition, the boundary of the flow set $C$ contains the jump set $D$ and its image $G(D)$, and the jump set $D$ is contained in a hyperplane, or a halfspace of a hyperplane. These assumptions are satisfied for a relevant class of hybrid systems, such as models of mechanical systems with impacts, see, for instance, the examples in Section~\ref{secexamples}.
\\
To be precise, we focus on the class of hybrid systems given by:
\begin{subequations}
\label{eqhsaf}
\bal{\label{eqhsafflow}
\dot x&= f(t,x),&x\in C,\\
x^+&=Lx+H,&x\in D \label{eqgd}
}
with the function $f$ measurable in its first argument and Lipschitz in its second argument, the matrix $L\in \mathbb{R}^{n\times n}$ being invertible, and $H\in \mathbb{R}^n$. Furthermore, the sets $C$ and $D$ are nonempty, closed and satisfy
\bal{\label{eqhsafc}
C&\subseteq \{x\in \mathbb{R}^n|\ Jx+K\leq 0\land \nonumber \\&\hspace{1.5cm}(JL^{-1} x+ K-JL^{-1}H)s \leq 0 \},\\
D&:=\{x\in C|\ Jx+K=0\land z_1x+z_2\leq 0\},\label{eqhsafd}
}
\end{subequations}
where the parameters $J^T,z_1^T\in \mathbb{R}^n\setminus \{0\}$, $K,z_2\in \mathbb{R}$ characterise the half hyperplane containing $D$,
and $s\in \{-1,1\}$ is selected such that $n_{gd}:=s(L^{-1})^TJ^T$ is a normal vector to $G(D)$ pointing out of $C$, cf.\ Fig.~\ref{figafps}, as we note that $G(D)\subset\{x\in \mathbb{R}^n|\ JL^{-1} x+  K-JL^{-1}H =0\}$ follows from the definitions of $D$ and $G$. Let $G(D)\subset C$ and the following assumption hold. 
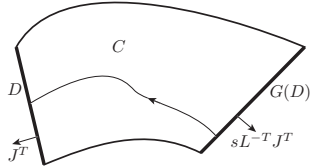
\fci{width=4cm}{sketch_linsys}{Pictorial illustration of the phase space of \eqref{eqhsaf} for the case $n=2$.\label{figafps}
}
\begin{ass}\label{assdgd}
The data of \eqref{eqhsaf} is such that there exist scalars $z_3,z_4,z_5>0$ such that
\begin{itemize}
\item $z_1x+z_2\geq z_3$ for all $x\in G(D)$,
\item $Jx+K<-z_4$ for all $x\in C$ that satisfy $|z_1x+z_2|\leq z_3$,
\item for all $x\in C$ with $z_1x+z_2\leq 0$, there exists a $y\in D$ such that $Jx+K\leq -z_5\|x-y\|$,
\item all maximal solutions of \eqref{eqhsaf} are complete.
\end{itemize}
%
%
\end{ass}
The first three bullets of this assumption are illustrated in Fig.~\ref{figsketch_ass3}.
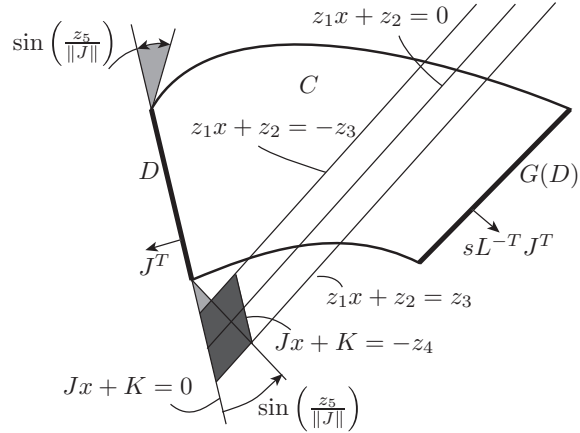
\fci{width=.5\columnwidth}{sketch_linsys_ass3}{Pictorial illustration of the phase space of \eqref{eqhsaf} when Assumption~\ref{assdgd} is satisfied. The second and third bullet of this assumption imply that the intersection between $C$ and the domains depicted in dark gray and light gray, respectively, is empty.\label{figsketch_ass3}
}
Note that this assumption directly implies $D\cap G(D)=\emptyset$, such that the first part of Assumption~\ref{asssysgtd} holds with $k=1$. In fact, $D$ and $G(D)$ are positioned at opposite sides of the hyperplane $\{x\in \mathbb{R}^n|\ z_1x+ z_2=0\}$.
We observe that all solutions to \eqref{eqhsaf} have a time domain that is unbounded in $t$-direction, as, firstly, $G(D)\cap D=\emptyset$ excludes Zeno-behaviour since $D$ is closed, secondly, $G$ is linear and, thirdly, $f$ is Lipschitz in its second argument. Hence, the hybrid system \eqref{eqhsaf} satisfies Assumptions~\ref{asssysgen} and \ref{asssysgtd}. In Section~\ref{secexamples}, we present examples of mechanical systems with impacts of the form \eqref{eqhsaf} that satisfy Assumption~\ref{assdgd}.
\\
In order to present a constructive Lyapunov function design, we first introduce the function $\bar G:\mathbb R^n\to \mathbb R^n$ as
 \beq{\label{eqbarg}
\bar G(x):=Lx\!+\!H\!+\!M(Jx\!+\!K)\!+\!sLJ^T\max(0,z_1x\!+\!z_2),
}
where the parameter $M\in \mathbb{R}^n$ is to be designed. Note that if $x\in D$, then $\bar G(x)=G(x)=Lx+H$.
\\
Since $G(D)\cap D=\emptyset$, Definition~\ref{dfnd} implies that $d(x,y)=0$ if and only if $x=y,$ or $ x=G(y)$ or $y=G(x)$. To design a Lyapunov function $V$, we note that \eqref{eqvbnds} requires that $V(x,y)=0$ if and only if $d(x,y)=0$. Hence, we propose the following piecewise quadratic Lyapunov function
\beq{
V(x,y)=
\min(\|x-y\|^2_{P_0},\|x-\bar G(y)\|^2_{P_s},\|\bar G(x)-y\|^2_{P_s}),
\label{eqv}
}
where the positive definite matrices $P_0,P_s\in \mathbb{R}^{n\times n}$ are to be designed. While this function is not smooth, we restrict our attention to a sufficiently small sub-level set where, as we will show in Lemma 4, the function is smooth.
\subsection*{Design of Lyapunov function parameters}
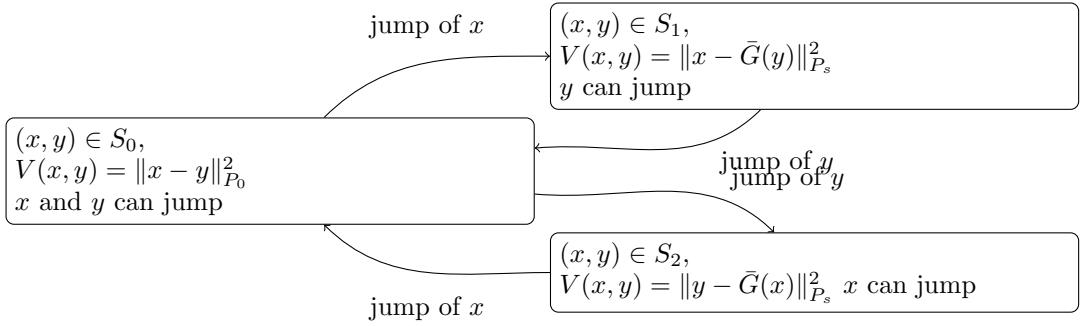
\begin{figure}
\begin{center}
\begin{tikzpicture}
\matrix[column sep=2mm,row sep=1mm]{
& \node[rounded corners=3pt,draw] (T){\begin{minipage}{.45\linewidth}
$(x,y)\in S_1$,\\
$V(x,y)=\|x-\bar G(y)\|_{P_s}^2$\\
$y$ can jump \end{minipage}};\\%
\node[rounded corners=3pt,draw] (L){\begin{minipage}{.45\linewidth}
$(x,y)\in S_0$,\\
$V(x,y)=\|x-y\|_{P_0}^2$\\
$x$ and $y$ can jump
 \end{minipage}};&
\\%
& \node[rounded corners=3pt,draw] (B){\begin{minipage}{.45\linewidth}
$(x,y)\in S_2$,\\
$V(x,y)=\|y-\bar G(x)\|_{P_s}^2$
$x$ can jump
\end{minipage}};\\%
}; 
\path[->](T) edge [out=-135,in=5] node[below right, near start](){jump of $y$} (L);
\path[->](L) edge[out=45,in=180] node[above=2mm](){jump of $x$} (T);
\path[<-](L) edge[out=-45,in=180] node[below=2mm](){jump of $x$} (B);
\path[->](L) edge [out=-5,in=135] node[above right , near end](){jump of $y$} (B);
\end{tikzpicture}
\caption{The three nodes indicate when $x$ and $y$ may jump provided $V(x,y)\leq v_L$, with $v_L$ as in Lemma~\ref{lemvcases}. When the conditions of Lemma~\ref{lemvdes} hold and, in addition, $V(x,y)\leq \max(1, e^{-\lambda d})v_L$ right before a jump, then this jump satisfies the scenarios depicted by arrows. \label{fig3j}
}
\end{center}
\end{figure}

To design the parameters $P_0,P_s$ and $M$ of the Lyapunov function $V$ in \eqref{eqv}, we distinguish the domains
\begin{subequations}\label{eqS012}
\bal{
S_0&:=\{(x,y)\in (C\!\cup\! D)^2|\ V(x,y)=\|x-y\|_{P_0}^2\},\\
S_1&:=\{(x,y)\in (C\!\cup\! D)^2|\ V(x,y)=\|x-\bar G(y)\|_{P_s}^2\},\\
S_2&:=\{(x,y)\in (C\!\cup\! D)^2|\ V(x,y)=\|\bar G(x)-y\|_{P_s}^2\}.
}\end{subequations}
The following lemma characterises the possibility of jumps from these sets, as illustrated in Fig.~\ref{fig3j}.
\begin{lem}\label{lemvcases}
Consider the hybrid system \eqref{eqhsaf} and $V$ in \eqref{eqv}. Let $M\in \mathbb{R}^n$ satisfy $(JL^{-1}M+1)s<0$, let  $P_0,P_s\succ 0$ and let Assumption~\ref{assdgd} hold.
There exists a $v_L>0$ such that for the sub-level set $V_L:=V^{-1}([0,v_L])$ the following statements hold:
\begin{itemize}
\item[1)]
$V_L$ consists of the union of the mutually disconnected sets $S_0\cap V_L$, $S_1\cap V_L$ and $S_2\cap V_L$.
\item[2)] $S_1\cap V_L\cap (D\times (C\cup D))=\emptyset$ and $S_2\cap V_L\cap (C\cup D)\times D)=\emptyset$.
\item[3)] it holds that $z_1y+z_2< 0$ for all $\bpm{x^T&y^T}^T\in S_1\cap V_L$ and $z_1x+z_2< 0$ for all $\bpm{x^T&y^T}^T\in S_2\cap V_L$.
\end{itemize}
\end{lem}
\begin{pf}
The proof is given in Appendix~\ref{apppfs}.
\end{pf}
We note that 3) implies that $\max(0,z_1x+z_2)=0$ in \eqref{eqbarg} if $V(x,y)=\|x-\bar G(y)\|_{P_s}^2<v_L$ or $V(x,y)=\|\bar G(x)-y\|_{P_s}^2<v_L$. In addition, this lemma allows to limit the number of jump scenarios between the sets $S_0,S_1,S_2$. Jumps of \eqref{eqhsaf} may trigger jumps between these cases. From item 2) in Lemma~\ref{lemvcases}, we observe that for $(x,y)\in S_1\cap V_L$ (or $(x,y)\in S_2\cap V_L$) jumps of $x$ (or $y$, respectively) are not feasible. In addition, when $x$ jumps and $(x,y)\in S_0$, then, after the jump, $G(x)\in G(D)$ implies $z_1G(x)+z_2\geq z_3$. When, in addition, $\|x-y\|_{P_0}^2$ was sufficiently small, then we find $\|G(x)-\bar G(y)\|_{P_s}^2<v_L$ (this statement has been made rigourous in the proof of the following lemma), and we obtain $(G(x),y)\in S_1$ after the jump. Consequently, the feasible jump scenarios are depicted in Fig.~\ref{fig3j}. Hence, \eqref{eqvjump} has to be proven along the four jumps depicted in Fig.~\ref{fig3j}.
For this reason, Lemma~\ref{lemvcases} is instrumental to prove the following lemma, that imposes the first design conditions on the parameters $P_0,P_s$ and $M$ of the Lyapunov function $V$ in \eqref{eqv}.
\begin{lem}\label{lemvdes}
Consider the hybrid system \eqref{eqhsaf}, let $M\in \mathbb{R}^n$ satisfy $(JL^{-1}M+1)s<0$, let  $P_0,P_s\succ 0$ and let Assumption~\ref{assdgd} hold. Consider the function $V$ in \eqref{eqv}.
If for some $\lambda_d\in \mathbb{R}$ it holds that
\bal{\label{lyapj1}
(L+MJ)^TP_s(L+MJ)\preceq e^{\lambda_d} P_0,\\
\label{lyapj2}
P_0\preceq e^{\lambda_d} P_s,
}
then there exist $\mathcal{K}_\infty$-functions $\alpha_1,\alpha_2$ and  $v_L>0$ such that the  conditions \eqref{eqvbnds} and \eqref{eqvjump} in Theorem~\ref{thmgen} are satisfied  with $V_L=V^{-1}([0,v_L])$ and
the function $V$ in \eqref{eqv} is smooth on an open domain containg $V_L$.
\end{lem}
\begin{pf}
The proof is given in Appendix~\ref{apppfs}.
\end{pf}

This lemma provides sufficient conditions on the dynamics of the hybrid systems such that the conditions on the Lyapunov function and its evolution along jumps of \eqref{eqhs} are satisfied. Additionally, \eqref{eqvflow} in Theorem~\ref{thmgen} imposes conditions on the evolution of the Lyapunov function along flows of \eqref{eqhs}. In the following section, we present a design for $V$ such that these conditions are also satisfied for an important class of hybrid systems.

\section{Tracking control problems\label{secapptrack}}
In this section, we will employ the results on the asymptotic stability of jumping hybrid trajectories to solve a tracking problem of a hybrid trajectory with jumps.
\\
We restrict our attention to tracking control problems for the class of systems \eqref{eqhsaf} with $f(t,x)=Ax+E+Bu(t,x)$, $A\in \mathbb{R}^{n\times n}$, $E,B\in \mathbb{R}^n$, with a control law $u: [0,\infty)\times C\to \mathbb{R}$ to be designed. In the scope of this tracking problem, we consider a reference trajectory $x_d$, which is a solution to \eqref{eqhsaf} for a feedforward input signal $u(x,t)=u_\text{ff}(t)$. We assume that $y$ is a trajectory that is generated by the control signal $u(t,y)=u_{\text{ff}}(t)+u_\text{fb}(t,y)$, and assume that $u_\text{fb}$ vanishes along the trajectory $x_d$, i.e.\, $u_\text{fb}(t,x_d(t))=0$ for almost all $t$. Hence, the flow map of the extended hybrid system \eqref{eqhse} is given by
\beq{\label{eqfeaf}
F_e(t,x_d,y)=\bpm{Ax_d\!+\!E\!+\!B(u_{\text{ff}}(t)\!+\!u_\text{fb}(t,x_d))\\Ay\!+\!E\!+\!B(u_{\text{ff}}(t)\!+\!u_\text{fb}(t,y))}.
}
Introducing the function $\displaystyle\bar x_d(t)\!:=\!x_d(t,\min_{(t,j)\in \dom x_d}j)$ design a switching feedback law $u_{fb}$ as:
\beq{\label{eqswitchcontrol}
u_\text{fb}(t,y)=\begin{cases}
-c_0(\bar x_d(t)-y),\\ \hspace{16mm} \mbox{for }\bpm{x_d^T(t)&y^T}^T\in S_0\\
-\frac{\beta_2^T}{\beta_2^T\beta_2}\beta_1(t)-c_1(\bar x_d(t)-\bar G (y)),\\ \hspace{16mm} \mbox{for }\bpm{x_d^T(t)&y^T}^T\in S_1\\
-\frac{\beta_4^T}{\beta_4^T\beta_4}\beta_3(t)-c_2(\bar G(\bar x_d(t))-y),\\
\hspace{16mm} \mbox{for }\bpm{x_d^T(t)&y^T}^T\in S_2
\end{cases}}
with $c_0^T,c_1^T,c_2^T\in \mathbb{R}^n$,
\bal{\nonumber
\beta_1(t)&=\bpm{I_n&-L-MJ}\bpm{A\bar x_d(t)+Bu_\text{ff}(t)+E\\A\bar G^{\circ}(\bar x_d(t))+Bu_\text{ff}(t)+E},
\\
\beta_3(t)&=\bpm{L+MJ& -I_n}\bpm{A\bar x_d(t)+Bu_\text{ff}(t)+E\\A\bar G(\bar x_d(t))+Bu_\text{ff}(t)+E}\nonumber, }
$\beta_2=-(L+MJ)B$ and $\beta_4=-B$, where $\bar G^\circ(x)$ is designed as $\bar G^\circ(x)=(L+MJ)^{-1}(x-H-MK)$, which coincides with the inverse of $\bar G$ restricted to $S_1\cap V_L$.
\\
Using this switched control law, which switches on the basis of the Lyapunov function designed in \eqref{eqv}, we formulate in the following result explicit conditions on the controller parameters $c_0,c_1,c_2, M,P_0$ and $P_s$ under which the tracking problem is solved.
\begin{thm}\label{thmtrack}
Consider the hybrid system \eqref{eqhsaf} with $f(t,x)=Ax+E+B(u_\text{ff}(t)+u_\text{fb}(t,x))$, for some measurable function $u_\text{ff}(t)$ and let $x_d$ be a solution of \eqref{eqhsaf} for $u_\text{fb}\equiv 0$.
Let $P_0,P_s\in \mathbb{R}^{n\times n}$, $M\in \mathbb{R}^n$, consider $V$ as in \eqref{eqv} and let $u_\text{fb}$ be designed as in \eqref{eqswitchcontrol}, with $\bar x_d(t)=x_d(t,\min_{(t,j)\in \dom x_d}j)$ and $c_0^T,c_1^T,c_2^T\in \mathbb{R}^{n}$. Let $L+MJ$ be invertible and $B\neq 0$.
\\
Let the assumptions of Lemma~\ref{lemvdes} hold for $\lambda_d\in \mathbb{R}$, let all trajectories of \eqref{eqhsaf} have a time domain that is unbounded in $t$-direction, and assume
\beq{\label{eqbetspans}
\beta_1(t)\in \mbox{span}(\beta_2),\quad\mbox{and } \beta_3(t)\in \mbox{span}(\beta_4)
}
hold for almost all $t$.
\\
Let, for some $\lambda_c\in \mathbb{R}$, the following LMIs be satisfied:
\bal{
(A+Bc_0)^TP_0+P_0(A+Bc_0)-\lambda_c P_0\preceq 0,\label{lmi1}\\
P_s(\beta_2c_1+(L+MJ)A(L+MJ)^{-1})+(\beta_2c_1+\nonumber\\
(L+MJ)A(L+MJ)^{-1})^TP_s+\lambda_cP_s\preceq 0\label{lmi3},\\
\label{lmi2}
P_s(A+Bc_2)+(A+Bc_2)^TP_s+\lambda_cP_s\preceq 0.
}
If either of the following cases hold, then the trajectory $x_d$ is asymptotically stable with respect to $d$.
\begin{enumerate}
 \item[1)] $\lambda_c<0,\lambda_d\leq 0$,
  \item[2)] all trajectories of \eqref{eqhs} have minimal average inter-jump time
$2\tau > 0$, $\lambda_c\leq 0$ and $\lambda_d+\lambda_c \tau<0$,
  \item[3)] all trajectories of \eqref{eqhs} have maximal average inter-jump time
$2\tau > 0$, $\lambda_d\leq 0$ and $\lambda_d+ \lambda_c \tau<0$.
\end{enumerate}
\end{thm}
\begin{pf}
The proof is given in Appendix~\ref{apppfs}.\end{pf}

\section{Examples \label{secexamples}}
We now present two examples to illustrate the results of this paper. In the first example, a tracking control problem is studied for a bouncing ball system where forces can be applied to the system. The tracking controller and Lyapunov function will be designed such that $\lambda_d=0$ and $\lambda_c<0$, such that case 1) of Theorem~\ref{thmtrack} is used to prove asymptotic stability of the reference trajectory.
\\
In the second example, a hybrid system is considered and a control law is proposed for which a maximal dwell-time argument proves asymptotic stability of the reference trajectory, illustrating case 3) of Theorem~\ref{thmtrack}.

\subsection*{Bouncing ball system with non-dissipative impacts}
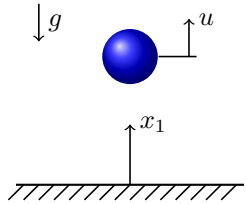
\begin{figure}[h!]
\begin{center}
\begin{tikzpicture}[-,auto,node distance=2cm,semithick]
\draw[thick] (-1.5,0)--(1.5,0);
\node[ball color=blue,circle,text=white] at (0,1.7) (A) {\phantom{M}};
\foreach \x in {-1.4,-1.2,...,1.4}
\draw (\x,0) -- (\x-.2,-.2);
\draw[-] (A.east)-- ++(0.5,0);
\draw[->] (A.east)+(0.4,0) -- ++ (0.4,0.5) node[right] {$u$};
\draw[->] (-1.2,2.4) --node[right]{$g$} (-1.2,1.9);
\draw[->] (0,0)--(0,.8) node[right]{$x_1$};
\end{tikzpicture}\hspace{1cm}
\caption{Bouncing ball system.\label{figbb}}
\end{center}
\end{figure}

We consider a tracking problem for a bouncing ball system with position $x_1$, unit mass, gravity $g=9.81$ and with impacts without energy dissipation, see Figure~\ref{figbb}.
\\
A finite control force $u$ can be applied to the ball along the vertical direction of motion and the constraint is positioned at $x_1=0$. The continuous-time dynamics can be described with

\beq{\ddot x_1=-g+u+\lambda,}
with $\lambda$ a finite constraint force satisfying the linear complementarity relation $\lambda\geq0, x_1\geq 0,\ \lambda x_1=0$. When this system arrives at $x_1=0$, then we apply Newton's impact law
\beq{
\dot x_1^+=-\dot x_1.
}
Provided that the finite constraint force $\lambda$  of this system can be ignored, (i.e., when for all time, either the position $x_1$ or the velocity $\dot x_1$ are non-zero) this system can be modelled as \eqref{eqhsaf} with $x=\bpm{x_1&\dot x_1}^T$,\\  $A=\bpm{0& 1\\0&0}$, $B=\bpm{0\\1}$, $E=\bpm{0\\-g}$, $L=-I_{2}$, $J=\bpm{1& 0}$, $H=0$, $K=0$, $s=-1$, $z_1=\bpm{0&1}$ and $z_2=0$ and the set $C$ is selected to exclude a neighbourhood of the origin, such that $D\cap G(D)=\emptyset$. While excluding a neighbourhood of the origin may imply that $C,D$ allow non-complete trajectories, in the following Lyapunov analysis, we prove invariance of a sub-level set $V_L=V^{-1}([0,v_L])$, for some $v_L$, of the Lyapunov function $V$ that characterises a neighbourhood of the reference trajectory and that does not contain points in $\partial C_e$ where solutions to \eqref{eqhse} cannot be extended. Hence, Assumptions~\ref{asssysgen}-\ref{assdgd} hold in the flow set $C_e\cap V_L$ and jump set $D_e\cap V_L$.
\\
We consider a tracking problem where the reference trajectory $x_d$ is a solution to \eqref{eqhsaf} for the feedforward function $u_{\text{ff}}(t)=0$ with initial condition $x_d(0,0)=\bpm{0& 10}^T$. Now, the tracking control law \eqref{eqswitchcontrol} is applied and Theorem~\ref{thmtrack} is used to find parameters $c_0,c_1,c_2,P_0,P_s,$ and $M$. We may select $\lambda_d=0$, $M=\bpm{0&0}^T$, $c_0=c_1=-\bpm{k&c}$, with $k=1,c=0.5$, and $P_0=P_s=\bpm{2.25&0.5\\0.5&2}$, that is the solution to $(A+Bc_0)^TP_0+P_0(A+Bc_0)=-I$. We select $\lambda_c=-0.25$, such that \eqref{lmi1} is satisfied, which, with $L+MJ=-I_2$, directly implies \eqref{lmi3}. Consequently,  Theorem \ref{thmtrack} proves that the trajectory $x_d$ is (locally) asymptotically stabilised with respect to $d$ by the control law \eqref{eqswitchcontrol}.
\\
In Fig.~\ref{figpeak}a-b), the reference trajectory $x_{d}$ and plant trajectory $x$ of the system is shown. In addition, the Euclidean error is depicted in panel c). Observe that $x$ accurately tracks the reference trajectory $x_d$ in the sense that jump times converge and, away from the impact times, the Euclidean error remains small and tends to zero. In Fig.~\ref{figuvper}, the applied control input $u$ and the Lyapunov function $V$ are shown.
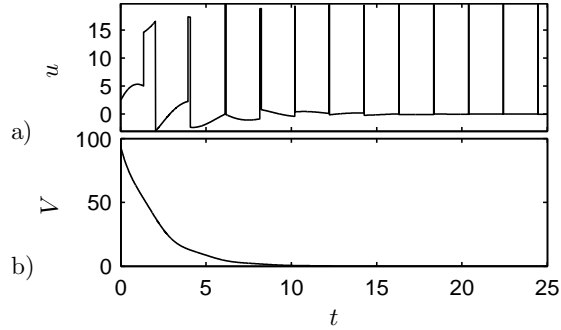
\fci{width=.5\columnwidth}{figuv}{a) Control action $u$ applied to the trajectory $x$ for the bouncing ball system with non-dissipative impacts. b) Lyapunov function $V$ in \eqref{eqv} evaluated along the trajectories shown in Fig.~\ref{figpeak}.\label{figuvper}}
In Fig.~\ref{figdnonpeak}, the distance function $d(x,x_d)$ in \eqref{eqd} is depicted when evaluated along these trajectories.
\\
Observe that this distance converges asymptotically to zero. Hence, the controller design in Section~\ref{secapptrack} solves the tracking problem and renders the trajectory $x_d$ asymptotically stable.
%
%
\begin{rmk}
In \cite{bie_wouw_13}, the same control law was designed in an ad-hoc manner for this specific case. Interestingly, the same control law $u_\text{fb}$ in \eqref{eqswitchcontrol} now follows in a systematic manner from the generic design framework presented in this paper. Clearly, this generic framework is applicable to a much wider range of examples, indicating the relevance of the presented work.
\end{rmk}

\subsection*{Dissipative mechanical system with impacts}
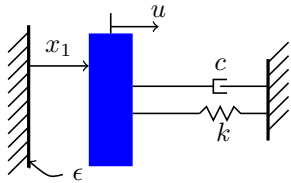
\begin{figure}[h!]
\begin{center}
\begin{tikzpicture}[-,auto,node distance=2cm,semithick,scale=.9]
\node[draw,color=blue,fill=blue,text=white] at (-1.8,0) (A) {\begin{minipage}{10mm} \phantom{M}\phantom{M}\phantom{M} \phantom{M}\end{minipage}};
\draw[-] (A.north)-- ++(0,0.3);
\draw[->] (A.north)+(0,.1) -- ++ (+.7,.1) node[above] {$u$};
\draw[-,very thick] (.5,-.6)--(.5,.6);
\draw[<-] (-3,-.9)..controls (-2.7,-1.15).. (-2.5,-1.1) node[right]{$\epsilon$};
\draw[<-] (A.west)+(0,.5)--node[above]{$x_1$} (-3,0.5);
\draw[-,very thick] (-3,-1)--(-3,1.1);
\foreach \y in {-.8,-.6,...,1.1}
\draw (-3,\y) -- (-3.3,\y-.3);
\foreach \y in {-.45,-.25,...,.65}
\draw (.5,\y) -- (.8,\y+.3);
\draw[-] (.5,.2)--(-.2,.2);
\draw[-] (-.1,.1)--(-.3,.1)--(-.3,.3)-- node[above]{$c$}(-.1,.3);
\draw[-] (A.east)+(0,.2)--(-.3,.2);
\draw[-] (-.5,-.2)--(-.4,-.1)--(-.3,-.3)--(-.2,-.1)--node[below]{$k$} (-.1,-.3) --(0,-.1)--(.1,-.2)--(.5,-.2);
\draw[-] (A.east)+(0,-.2)--(-.5,-.2);
\end{tikzpicture}\hspace{1cm}
\caption{Dissipative mechanical system.\label{figsys2}}
\end{center}
\end{figure}

Now, we consider a single degree-of-freedom system with a damper with  damping constant $c>0$ and a spring with stiffness $k>0$ and unloaded position $x=\bar x_1$, as shown in Fig.~\ref{figsys2}. Impacts can only occur at the constraint at $x_1=0$.
Let the impacts be described by a restitution coefficient $\varepsilon=0.9$. Hence, the impacts are dissipative, which allows to study the stability of the trajectory using a maximal average inter-jump time result. Assuming that finite constraint forces can be ignored, the hybrid system is described by \eqref{eqhsaf} with\\ $A=\bpm{0& 1\\ -k&-c}$, $B=\bpm{0\\1}$, $E=\bpm{0\\k\bar x_1}$, $L=-\varepsilon I_2$, $J=\bpm{1&0}$, $K=0$, $H=0$, $s=-1$, $z_1=\bpm{0&1}$, $z_2=0$  and the set $C$ is selected to exclude the origin. The parameters $\bar x_1=1$,  $k=1$ and $c=0.02$ are used.
\\
Let the reference trajectory $x_d$ be a solution to \eqref{eqhsaf} for a feedforward function $u=u_{\text{ff}}(t)=100\cos(\omega t)$, with $\omega=0.4$, as shown in Fig.~\ref{figepsnc}. This forcing is selected such that the reference trajectory $x_d$ with initial condition $x_d(0,0)=\bpm{50&0}^T$  has a maximal average inter-jump time $\tau_d>0$. In addition, $x_d(t,j)$  does not arrive at the origin. To show that the tracking problem is not trivial, in the same figure, a trajectory of \eqref{eqhsaf} is depicted with the same forcing and initial condition $x_d(0,0)=\bpm{51&0}^T$. Clearly, both trajectories diverge.
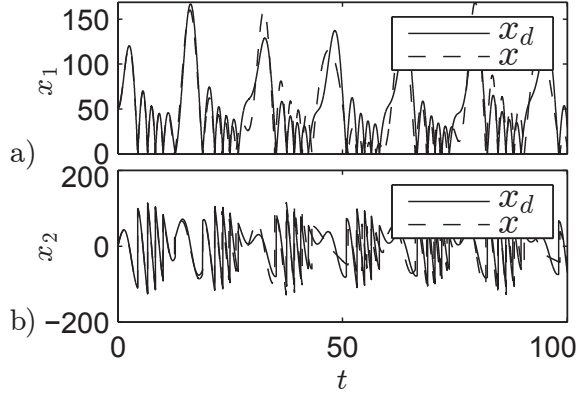
\fci{width=.5\columnwidth}{figepsnc}{Reference trajectory $x_d$ for the dissipative mechanical system and trajectory $x$ with nearby initial condition for $u=u_{\text{ff}}(t)$.\label{figepsnc}}
\\
For this example, we derive an explicit expression for the distance function. As the reference trajectory $x_d$ stays away from the origin, we can model the system with $D=\{0\}\times (-\infty,-r]$, with $r>0$ sufficiently small, such that $G(D)=\{0\}\times [\varepsilon r,\infty)$. Then, \eqref{eqd} yields
\beq{\label{eqd123}
d(x,y)=\min(d_0(x,y),d_1(x,y),d_2(x,y)),
}
with
$d_0(x,y)=\inf_{z_x=z_y\in C}\left\|\bpm{x-z_x\\y-z_y}\right\|$,
$d_1(x,y)=\inf_{z_x=G(z_y), z_y\in D}\left\|\bpm{x-z_x\\y-z_y}\right\|$ and $d_2(x,y)=d_1(y,x)$. We will now find explicit expressions for $d_0$ and $d_1$. By the observation that $d_0^2(x,y)=\inf_{z_x\in C}\|x-z_x\|^2+\|y-z_x\|^2$, we observe that the infimum is attained $z_x=\frac{x+y}2$, which is contained in $C$ when, firstly, $d_0(x,y)$ is sufficiently small and secondly, $x=x_d$ with the reference trajectory sufficiently far from the origin. Hence, we obtain
\bal{
d_0(x,y)=\frac{1}{\sqrt 2}\|x-y\|\label{eqd0}.
}
Now, we compute $d_1(x,y)$ using that
\bal{\label{eqd1eps}
d_1(x,y)^2&=\inf_{z_x=G(z_y), z_y\in D}\left\|\bpm{x-z_x\\y-z_y}\right\|^2\\
&=\inf_{z_{y2}\leq -r}x_1^2+y_1^2+(x_2+\varepsilon z_{y2})^2+(y_2-z_{y2})^2.\nonumber
}
Studying $\frac{d}{dz_{y2}}(x_1^2+y_1^2+(x_2+\varepsilon z_{y2})^2+(y_2-z_{y2})^2)=2\varepsilon(x_2+\varepsilon z_{y2})-2(y_2-z_{y2})=0$, leads to the minimiser $z_{y2}=\frac{y_2-\varepsilon x_2}{1+\varepsilon^2}$ when $\frac{y_2-\varepsilon x_2}{1+\varepsilon^2}<-r$, and $z_{y_2}=-r$ otherwise. Substituting this expression in \eqref{eqd1eps}, we obtain, after some algebraic manipulations, that
\bal{
d_1(x,y)&=
\begin{cases}
\sqrt{x_1^2+y_1^2+\frac{(\varepsilon y_2+x_2)^2}{1+\varepsilon^2}},
\\\hspace{1cm}\mbox{ if }\frac{y_2-\varepsilon x_2}{1+\varepsilon^2}<-r,\\
\sqrt{x_1^2\!+\!y_1^2\!+\!(x_2\!-\!\varepsilon r)^2\!+\!(y_2\!+\!r)^2},
\\\hspace{1cm}\mbox{otherwise,}
\label{eqd1expl}
\end{cases}
}
where, in the first case, we used the relations $x_2+\varepsilon z_{y2}=\frac{x_2+\varepsilon y_2}{1+\varepsilon^2}$ and $y_2-z_2=\frac{\varepsilon (x_2+\varepsilon y_2)}{1+\varepsilon^2}$. We note that for $\varepsilon>0$, $r\ll \|x\|,\|y\|$ implies that the second case in \eqref{eqd1expl} will not be attained when $d(x,y)=d_1(x,y)$. Recall that $d_1(x,y)=d_2(y,x)$, such that \eqref{eqd123} gives \beq{
d(x,y)=\min(d_0(x,y),d_1(x,y),d_1(y,x)),
}
with $d_0$ in \eqref{eqd0} and $d_1$ in \eqref{eqd1expl}. 
\\
We now apply the constructive control law design proposed in Section~\ref{secapptrack} to enforce tracking of the trajectory $x_d$. Selecting $P_0=\bpm{k&0\\0&1}$ and $P_s=\tfrac 1{\varepsilon}P_0$, we observe that the conditions of Lemma~\ref{lemvdes} are satisfied with $\lambda_d=\log(\varepsilon)<0$. In addition, we observe that $c_0=c_1=c_2=0$ can be selected, such that \eqref{lmi1}-\eqref{lmi2} hold with $\lambda_c=0$, as $P_0A+A^TP_0=\bpm{0& 0\\0& -2c}$ and $P_sA+A^TP_s=\bpm{0& 0\\0& -\frac{2c}{\varepsilon^2}}$. Then, \eqref{eqswitchcontrol} yields the control law:
\beq{\label{eqswitchred}
u_\text{fb}(t,x)=\begin{cases}
0,\quad  \bpm{\bar x_d^T(t)&y^T}^T\in S_0\\
-\frac{1+\varepsilon}{\varepsilon}(k\bar x_1+u_{\text{ff}}(t)),\\ \hspace{1cm} \bpm{\bar x_d^T(t)&y^T}^T\in S_1\\
-(1+\varepsilon)(k\bar x_1+u_{\text{ff}}(t)),\\ \hspace{1cm} \bpm{\bar x_d^T(t)&y^T}^T\in S_2.
\end{cases}}
As the trajectory $x_d$ has a maximal average inter-jump time, denoted $\tau_d$, nearby trajectories will have the same behaviour. Hence, selecting $v_L>0$ sufficiently small and restricting our attention to the hybrid system \eqref{eqhse} with flow set $C_e\cap V_L$ and jump set $D_e\cap V_L$, with $V_L=V^{-1}([0,v_L])$, we conclude that $x$ also has a maximal average dwell-time $\tau_x$, with $\tau_x$ close to $\tau_d$. Hence, the trajectory of the embedded system \eqref{eqhse} has a maximal average inter-jump time $\frac{\max(\tau_d,\tau_x)}{2}>0$. Consequently, case 3) of Theorem~\ref{thmtrack} proves that the trajectory is (locally) asymptotically stabilised with respect to $d$ by the control law \eqref{eqswitchred}.
\\
In Fig.~\ref{figeps}, the performance of this controller is illustrated and  a trajectory with initial condition $x(0,0)=\bpm{100&0}^T$ is shown to $x_d$.
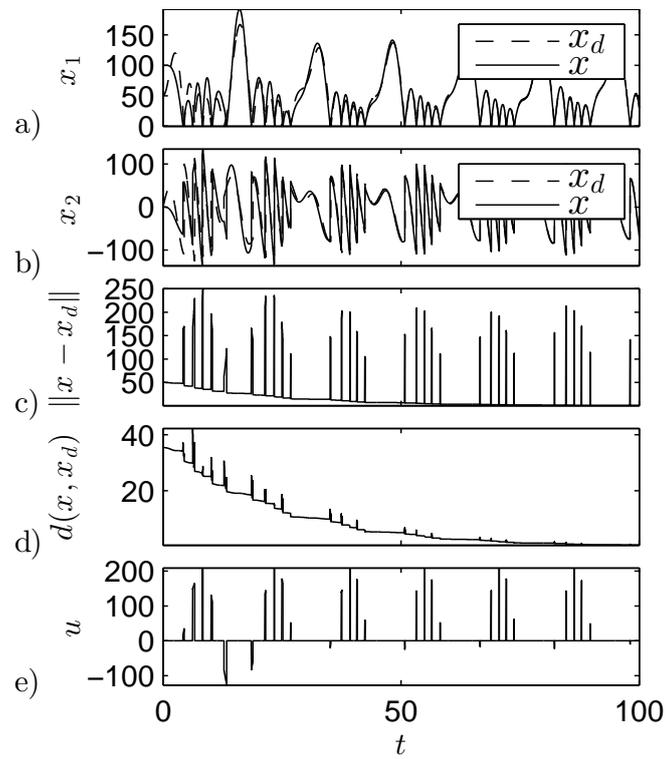
\fci{width=.6\columnwidth}{figeps}{a) and b) Reference trajectory $x_d$ and plant trajectory $x$ for the dissipative mechanical system and periodic forcing. c)~Euclidean tracking error. d)~Distance function \eqref{eqd}. e)~Control force $u$.\label{figeps}}
\\
From the structure of the control law \eqref{eqswitchred}, we observe that no control is active when $V(\bar \varphi_y(t,j),x_d(t,j))=\|\varphi_y(t,j)-x_d(t,j)\|_{P_0}$. In fact, the dissipative effect of both the damping force $c\dot x$ and the jump map implies that no control is needed during these time intervals. The control input $u$ only needs to compensate the potentially destabilising effect of the forcing term $(E+Bu_{\text{ff}})$ during the ``peaks'' of the Euclidean error.

\section{Conclusion\label{secconcl}}
In this paper, we considered the stability of time-varying and jumping trajectories of hybrid systems with state-triggered jumps, which is essential in tracking control, observer design and synchronisation problems.
A general distance function design was proposed that allows to compare two trajectories of a hybrid system, thereby enabling the stability analysis for hybrid trajectories. This stability analysis was employed and allowed to propose a constructive tracking control law for a class of hybrid systems.
%
%
\\
The proposed design for the distance function takes the nature of the jumps of the hybrid system into account, such that the distance function provides a good comparison between the trajectories of two solutions, without the ``peaking behaviour'' along solutions that is expected in the Euclidean distance. The stability properties of trajectories were studied in terms of this distance function. Sufficient conditions for stability have been formulated using Lyapunov functions with sub-level sets that consist of disconnected pieces. Moreover, the conditions are formulated in terms of  maximum or minimum average inter-jump time conditions to allow for increase of the Lyapunov function over flow or jumps, respectively.
In case the jump map is an affine function and the jump set a hyperplane, a piecewise quadratic Lyapunov function was proposed that can be constructed systematically.
\\
Focussing on a class of tracking problems for hybrid systems where control is only possible during flow, we designed a switching tracking control law. The asymptotic stability of a reference trajectory of the closed-loop system can be assessed with matrix conditions obtained from our general theory.
\\
Finally, we applied the proposed tracking control law in two examples and observed that the control law achieves accurate tracking. These examples also illustrated that the presented asymptotic stability notion does correspond to desired tracking behaviour. This underlines that the proposed distance function enables a good comparison between hybrid trajectories and has the potential to play an important role in tracking control, observer design and synchronisation problems for hybrid systems.
\subsection*{Acknowledgement}
J.J.B.~Biemond is a FWO Pegasus Marie Curie Fellow. This research is supported in part by the European Union Seventh Framework Programme [FP7/2007-2013]  under grant agreement no. 257462 HYCON2 Network of excellence.

\bibliographystyle{plain}
\bibliography{C:/kuleuven/literature/allrefs}

\providecommand{\noopsort}[1]{}
\begin{thebibliography}{10}

\bibitem{ang_02}
D.~Angeli.
\newblock A {Lyapunov} approach to incremental stability properties.
\newblock {\em IEEE Transactions on Automatic Control}, 47(3):410--421, 2002.

\bibitem{aub_fra_09}
J.-P. Aubin and H.~Frankowska.
\newblock {\em Set-Valued Analysis}.
\newblock Modern Birkh\"auser Classics. Birkh\"auser, Boston, 2009.

\bibitem{bie_wouw_13}
J.~J.~B. Biemond, N.~van~de Wouw, W.~P. M.~H. Heemels, and H.~Nijmeijer.
\newblock Tracking control for hybrid systems with state-triggered jumps.
\newblock {\em IEEE Transactions on Automatic Control}, 58(4):876--890, 2013.

\bibitem{bro_99}
B.~Brogliato.
\newblock {\em Nonsmooth Mechanics}.
\newblock Springer-Verlag, London, 1999.

\bibitem{bro_ara_02}
M.~Broucke and A.~Arapostathis.
\newblock Continuous selections of trajectories of hybrid systems.
\newblock {\em Systems \and Control Letters}, 47:149--157, 2002.

\bibitem{for_teel_11tac}
F.~Forni, A.~R. Teel, and L.~Zaccarian.
\newblock Follow the bouncing ball: global results on tracking and state
  estimation with impacts.
\newblock {\em submitted to IEEE Transactions Automatic Control}, 2011.

\bibitem{for_teel_11_1}
F.~Forni, A.~R. Teel, and L.~Zaccarian.
\newblock Tracking control in billiards using mirrors without smoke, {Part}
  {I}: {Lyapunov}-based local tracking in polyhedral regions.
\newblock In {\em Proceedings of the 50th IEEE Conference on Decision and
  Control, Orlando}, pages 3283--3288, 2011.

\bibitem{for_teel_11_2}
F.~Forni, A.~R. Teel, and L.~Zaccarian.
\newblock Tracking control in billiards using mirrors without smoke, {Part}
  {II}: additional {Lyapunov-based} local and global results.
\newblock In {\em Proceedings of the 50th IEEE Conference on Decision and
  Control, Orlando}, pages 3289--3294, 2011.

\bibitem{gal_men_12}
S.~Galeani, L.~Menini, and A.~Potini.
\newblock Robust trajectory tracking for a class of hybrid systems: an internal
  model principle approach.
\newblock {\em IEEE Transactions on Automatic Control}, 57(2):344--359, 2012.

\bibitem{gal_men_08}
S.~Galeani, L.~Menini, A.~Potini, and A.~Tornamb\`e.
\newblock Trajectory tracking for a particle in elliptical billiards.
\newblock {\em International Journal of Control}, 81(2):189--213, 2008.

\bibitem{goe_san_12}
R.~Goebel, R.~G. Sanfelice, and A.~R. Teel.
\newblock {\em Hybrid dynamical systems: Modeling, Stability and Robustness}.
\newblock Princeton University Press, Princeton, 2012.

\bibitem{hee_cam_11}
W.~P. M.~H. Heemels, M.~K. Camlibel, J.~M. Schumacher, and B.~Brogliato.
\newblock Observer-based control of linear complementarity systems.
\newblock {\em International Journal of Robust and Nonlinear Control},
  21(10):1193--1218, 2011.

\bibitem{hee_sch_10}
W.~P. M.~H. Heemels, B.~De~Schutter, J.~Lunze, and M.~Lazar.
\newblock Stability analysis and controller synthesis for hybrid dynamical
  systems.
\newblock {\em Philosophical Transactions of the Royal Society A:
  Mathematical,Physical and Engineering Sciences}, 368(1930):4937--4960, 2010.

\bibitem{hes_lib_08}
J.~P. Hespanha, Daniel Liberzon, and Andrew~R. Teel.
\newblock {Lyapunov} conditions for input-to-state stability of impulsive
  systems.
\newblock {\em Automatica}, 44(11):2735 -- 2744, 2008.

\bibitem{kha_02}
H.~K. Khalil.
\newblock {\em Nonlinear Systems}.
\newblock Prentice Hall, Upper Saddle River, third edition, 2002.

\bibitem{kol_fom_75}
A.~N. Kolmogorov and S.~V. Fomin.
\newblock {\em Introductory real analysis}.
\newblock Dover, New York, 1975.

\bibitem{lei_wouw_08}
R.~I. Leine and {\noopsort{Wouw}}{N. van de Wouw}.
\newblock {\em Stability and convergence of mechanical systems with unilateral
  constraints}, volume~36 of {\em Lecture Notes in Applied and Computational
  Mechanics}.
\newblock Springer-Verlag, Berlin, 2008.

\bibitem{lei_wouw_08ijbc}
R.~I. Leine and {\noopsort{Wouw}}{N. van de Wouw}.
\newblock Uniform convergence of monotone measure differential inclusions: with
  application to the control of mechanical systems with unilateral constraints.
\newblock {\em International Journal of Bifurcation and Chaos},
  18(5):1435--1457, 2008.

\bibitem{leo_08}
G.~A. Leonov.
\newblock {\em Strange attractors and classical stability theory}.
\newblock St.~Petersburg University Press, St.~Petersburg, 2008.

\bibitem{li_phi_14}
Y.~Li, S.~Phillips, and R.~G. Sanfelice.
\newblock Results on incremental stability for a class of hybrid systems.
\newblock In {\em to appear in Proceedings of the 2014 IEEE Conference on
  Decision and Control, Los Angelos}, 2014.

\bibitem{lun_lam_09}
J.~Lunze and F.~Lamnabhi-Lagarrigue, editors.
\newblock {\em Handbook of hybrid systems control}.
\newblock Cambridge University Press, Cambridge, 2009.

\bibitem{lyg_joh_03}
J.~Lygeros, K.~H. Johansson, S.~N. Simi\'c, Jun Zhang, and S.~S. Sastry.
\newblock Dynamical properties of hybrid automata.
\newblock {\em IEEE Transactions on Automatic Control}, 48(1):2--17, 2003.

\bibitem{men_tor_01}
L.~Menini and A.~Tornamb\`e.
\newblock Asymptotic tracking of periodic trajectories for a simple mechanical
  system subject to nonsmooth impacts.
\newblock {\em IEEE Transactions on Automatic Control}, 46(7):1122--1126, 2001.

\bibitem{mor_78}
J.~J. Moreau.
\newblock Approximation en graphe d'une \'evolution discontinue.
\newblock {\em RAIRO Analyse Num\'erique}, 12(1):75--84, 1978.

\bibitem{mor_bro_10}
I.~C. Mor\u{a}rescu and B.~Brogliato.
\newblock Trajectory tracking control of multiconstraint complementarity
  {Lagrangian} systems.
\newblock {\em IEEE Transactions on Automatic Control}, 55(6):1300--1313, 2010.

\bibitem{pos_tab_14}
R~Postoyan, P~Tabuada, D~Nesic, and A~Anta.
\newblock A framework for the event-triggered stabilization of nonlinear
  systems.
\newblock {\em IEEE Transactions on Automatic control}, to appear, 2014.

\bibitem{ruf_wouw_13}
B.~S. R\"uffer, N.~van~de Wouw, and M.~Mueller.
\newblock Convergent systems vs. incremental stability.
\newblock {\em Systems \& Control Letters}, 62(3):277--285, 2013.

\bibitem{san_bie_14}
R.~G. Sanfelice, J.~J.~B. Biemond, N.~van~de Wouw, and W.~P. M.~H. Heemels.
\newblock An embedding approach for the design of state-feedback tracking
  controllers for references with jumps.
\newblock {\em International Journal of Robust and Nonlinear Control},
  24(11):1585--1608, 2014.

\bibitem{san_goe_07}
R.~G. Sanfelice, R.~Goebel, and A.~R. Teel.
\newblock Invariance principles for hybrid systems with connections to
  detectability and asymptotic stability.
\newblock {\em IEEE Transactions on Automatic Control}, 52(12):2282--2297,
  2007.

\bibitem{san_goe_08}
R.~G. Sanfelice, R.~Goebel, and A.R. Teel.
\newblock Generalized solutions to hybrid dynamical systems.
\newblock {\em ESAIM: Control, Optimisation and Calculus of Variations},
  14(4):699--724, 2008.

\bibitem{sch_sch_00}
{\noopsort{Schaft}}{A. J. van der Schaft} and J.~M. Schumacher.
\newblock {\em An introduction to hybrid dynamical systems}, volume 251 of {\em
  Lecture Notes in Control and Information Sciences}.
\newblock Springer-Verlag, London, 2000.

\bibitem{sch_98}
M.~Schatzman.
\newblock Uniqueness and continuous dependence on data for one-dimensional
  impact problems.
\newblock {\em Mathematical and Computer Modelling}, 28(4-8):1--18, 1998.

\bibitem{zam_wouw_13}
M.~Zamani, N.~van~de Wouw, and R.~Majumdar.
\newblock Backstepping controller synthesis and characterizations of
  incremental stability.
\newblock {\em Systems \& Control Letters}, 62(10):949--962, 2013.

\end{thebibliography}

\appendix

\section{Alternative distance function
\label{secaltdist}}
The distance function \eqref{eqd} is not necessarily continuous over jumps, when evaluated along solutions to \eqref{eqhs}. When $G$ is a single-valued and invertible function, such a continuity property could be induced by the alternative distance function design:
\beq{\label{eqdq}
d_{Q}(x,y)=\inf_{N\in \mathbb{N}} \inf_{\substack{ \bpm{x^{iT},y^{iT}}^T\in \mathcal{A},\ i=1,\ldots, N,\\y^0=x,\ x^{N+1}=y}}\sum_{i=0}^N\|y^i-x^{i+1}\|,
}
that coincides with the quotient metric, cf.\ \cite{kol_fom_75}, on the quotient space generated by the equivalence $x\sim y$ if $(x,y)\in \mathcal{A}$. This quotient space has been suggested in \cite{lyg_joh_03} to study hybrid systems. We note that when $G$ is non-invertible, then $d_Q(x,y)=0 \Leftrightarrow \bpm{x^T&y^T}^T\in \mathcal{A}$ may not hold. To allow for non-invertible jump maps, we prefer the distance function $d$ in \eqref{eqd} over $d_Q$ in \eqref{eqdq}.

\section{Proofs\label{apppfs}}
To prove Lemma~\ref{lemvcases} we need the following result.
\begin{lem}\label{lembarg}
Consider the hybrid system \eqref{eqhsaf} and $\bar G$ in \eqref{eqbarg}, let $M\in \mathbb{R}^n$ satisfy $(JL^{-1}M+1)s<0$, let  $P_0,P_s\succ 0$ and let Assumption~\ref{assdgd} hold.
There exists a $\delta_1>0$ such that if $\|x-\bar G(y)\|\leq \delta_1$ and $x,y\in C\cup D$, then $z_1y+z_2\leq -z_3$ and $z_1\bar G(y)+z_2>\frac{z_3}2$.
\end{lem}
\begin{pf}
To prove the lemma, we first present a condition on $\delta_1$, such that $\|x-\bar G(y)\|> \delta_1$ holds for all $x,y\in C\cup D$ with $z_1y+z_2> -z_3$, thereby proving the first implication in the lemma.
\\
We introduce the intermediate variables $z=y+\tfrac{K J^T}{JJ^T}$, $\mu=\max(0,z_1y+z_2)\geq 0$, the vector $n_{gd}=s(L^{-1})^TJ^T$, scalar $k_{gd}=sK-sJL^{-1}H$ and affine function
\[
c_{gd}(x)=n_{gd}^Tx+k_{gd}=s(JL^{-1}x+K-JL^{-1}H),
\]
 such that for all $x\in C$, $c_{gd}(x)\leq 0$, cf.\ \eqref{eqhsafc}.
To provide a lower bound on $\|x-\bar G(y)\|$, we will first derive a lower bound on $c_{gd}(\bar G(y))$. Observe that
\bal{
c_{gd}(\bar G(y))=&s(JL^{-1}\bar G(y)+K-JL^{-1}H)\nonumber \\
=&s\bigg(JL^{-1}\bigg( (L+MJ)\left(\frac{-K J^T}{JJ^T}+z\right)\nonumber \\&+MK+sLJ^T\mu\bigg)+K\bigg)\nonumber \\
=&s
(1+JL^{-1}M)Jz
+s^2JJ^T\mu\label{eqcgdjz}
}
with $\mu=\max(0,z_1y+z_2)\geq 0$. We directly observe $c_{gd}(\bar G(y))\geq 0$ as
$s(1+JL^{-1}M)<0$ holds by the assumption in the lemma, and $Jz=Jy+K\leq 0$ for $y\in C$.
\\
We will now provide a positive lower bound on $c_{gd}(\bar G(y))$ for $y\in C$ with $z_1y+z_2> -z_3$.
By Assumption~\ref{assdgd},
\beq{\label{eqz4}
Jz=Jy+K<-z_4\quad \mbox{for } y\in C\mbox{ satisfying } |z_1y+z_2|\leq z_3,}
holds, with $z_4>0$ given in the assumption.
From  \eqref{eqcgdjz} and \eqref{eqz4}, we conclude that $c_{gd}(\bar G(y))>-s
(1+JL^{-1}M)z_4>0$ if $|z_1y+z_2|\leq z_3$ and $y\in C$. When $z_1y+z_2\geq z_3$, we obtain $\mu\geq z_3$ and $c_{gd}(\bar G(y))\geq JJ^Tz_3>0$. Hence, we have shown
\bal{
c_{gd}(\bar G(y))\geq& \min(-s(1+JL^{-1}M)z_4, JJ^Tz_3),\nonumber \\
&\mbox{ for }y\in C \mbox{ with }z_1y+z_2> -z_3.\label{eqz3}
}
For all $x,y\in C$ it follows from the Cauchy-Schwarz inequality that
\bal{
\|x-\bar G(y)\|
&\geq \tfrac{1}{\|n_{gd}\|}| n_{gd}^T(x-\bar G(y))|,\nonumber \\
&= \tfrac{1}{\|n_{gd}\|}|c_{gd}(x)-c_{gd}(\bar G(y))|,\nonumber \\
\intertext{
and, as  $c_{gd}(\bar G(y))\geq 0$ and $c_{gd}(x)<0$, we find
}
\|x-\bar G(y)\|
&\geq \tfrac{1}{\|n_{gd}\|}c_{gd}(\bar G(y)),\mbox{ for }x,y\in C,\label{eqncxg}\\
\intertext{and, using \eqref{eqz3}, }
\|x-\bar G(y)\|&\geq \tfrac{1}{\|JL^{-1}\|} \min (-s
(1+JL^{-1}M)z_4,JJ^Tz_3),\nonumber \\
&\hspace{15mm} \mbox{for }x,y\in C,\ z_1y+z_2\geq -z_3.\nonumber
}
Hence, if $x,y\in C$ and $\|x-\bar G(y)\|\leq \delta_1$ with
\beq{
\delta_1<\tfrac{1}{\|JL^{-1}\|} \min (-s
(1+JL^{-1}M)z_4,JJ^Tz_3).\label{conddel1a}
}
then $z_1y+z_2<-z_3$ holds.
\\
We now derive an upper bound on $\delta_1$ satisfying \eqref{conddel1a} such that for $x,y\in C$ with $\|x-\bar G(y)\|\leq \delta_1$, the relation $z_1\bar G(y)+z_2\geq \tfrac{z_3}{2}$ holds. As shown above, $z_1y+z_2<-z_3$ holds, such that Assumption~\ref{assdgd} implies that there exists $y^\star\in D$ such that $\|y-y^\star\|\leq \tfrac{Jy+K}{-z_5}$. With \eqref{eqcgdjz}, we find $Jy+K=Jz=\tfrac{c_{gd}(\bar G(y))}{s(1+JL^{-1}M)}$ (as $\mu=0$), such that
\beq{
\|y-y^\star\|\leq
\tfrac{c_{gd}(\bar G(y))}{-z_5s(1+JL^{-1}M)}
}
is obtained. Using \eqref{eqncxg}, we observe that
\[
\|y-y^\star\|\leq
\tfrac{\|n_{gd}\| \|x-\bar G(y)\|}{-z_5s(1+JL^{-1}M)}\leq
\tfrac{\|n_{gd}\|\delta_1}{-z_5s(1+JL^{-1}M)}
\]
holds for $\|x-\bar G(y)\|\leq \delta_1$.
Denoting the largest singular value of $(L+MJ)$ with $\ell_g$, we observe that $\ell_g>0$ since, otherwise, $L+MJ=0$, such that $\frac{s}{JJ^T}JL^{-1}(L+MJ)J^T=s(1+JL^{-1}M)=0$, which contradicts the assumption in the lemma. We find
\[
\|\bar G(y)-\bar G(y^\star)\|
\leq \tfrac{l_g \|n_{gd}\|\delta_1}{-z_5s(1+JL^{-1}M)}.
\]
With the relation
\[
z_1\bar G(y)+z_2\geq z_1\bar G(y^\star)+z_2-\|z_1\| \|\bar G(y)-\bar G(y^\star)\|
\]
and the observation that $y^\star\in D$ implies $z_1\bar G(y^\star)+z_2\geq z_3$ by Assumption~\ref{assdgd}, we find
\[
z_1\bar G(y)+z_2\geq z_3-\|z_1\| \tfrac{l_g \|n_{gd}\|\delta_1}{-z_5s(1+JL^{-1}M)}
\]
such that $z_1\bar G(y)+z_2\geq \frac{z_3}{2}$ if
\beq{
\delta_1\leq \tfrac{-z_5s(1-JL^{-1}M)z_3}{2\|z_1\| l_g \|n_{gd}\|},\label{conddel1b}
}
whose right-hand side is positive as $s(1-JL^{-1}M)<0$.
\\
Hence, for any $\delta_1>0$ satisfying \eqref{conddel1a} and \eqref{conddel1b},  $x,y\in C$ and $\|x-\bar G(y)\|\leq \delta_1$ implies $z_1y+z_2<-z_3$ and $z_1\bar G(y)+z_2\geq \frac{z_3}{2}$, thereby proving the lemma.
\end{pf}

\begin{pf}[Proof of Lemma~\ref{lemvcases}]
In order to prove the first part of the lemma, we select a $\delta_1<\frac{z_3}{2\|z_1\|}$ satisfying the conditions of Lemma~\ref{lembarg}, and fix $v_L>0$ satisfying $v_L<{\underline \lambda}{ \delta_1^2}$ and
$v_L<{\underline \lambda}{(\tfrac{3z_3}{2\|z_1\|}-\delta_1)^2}$, with $\underline \lambda>0$ smaller than the eigenvalues of $P_0$ and $P_s$.
\\
We now consider points $x,y\in C$ such that $\|x-\bar G(y)\|_{P_s}^2\leq v_L$ and prove that this implies $\|x-y\|_{P_0}^2>v_L$ and $\|\bar G(x)-y\|_{P_s}^2>v_L$. By the definition of $\underline \lambda$, we observe that $\|x-\bar G(y)\|_{P_s}^2\leq v_L$ implies $\|x-\bar G(y)\|<\delta_1$, such that Lemma~\ref{lembarg} implies $z_1y+z_2\leq -z_3$ and $z_1\bar G(y)+z_2>\frac{z_3}{2}$. Hence, $\|y-\bar G(y)\|\geq \frac{|z_1 (y-\bar G(y))|}{\|z_1\|}\geq \frac{3z_3}{2\|z_1\|}$ holds. We then obtain
\bal{
\|x-y\|_{P_0}^2&\geq \underline \lambda \|x-y\|^2\geq \underline \lambda (\|y-\bar G(y)\|-\|x-\bar G(y)\|)^2\nonumber \\
&\geq \underline \lambda (\tfrac{3z_3}{2\|z_1\|}-\delta_1)^2>v_L.
}
To prove $\|\bar G(x)-y\|_{P_s}^2>v_L$ when $\|x-\bar G(y)\|_{P_s}^2\leq v_L$, we suppose, for the sake of contradiction, that $\|\bar G(x)-y\|_{P_s}^2\leq v_L$ and $\|x-\bar G(y)\|_{P_s}^2\leq v_L$. Then, following analogous reasoning as above, from $\|\bar G(x)-y\|_{P_s}^2\leq v_L$ we conclude $z_1x+z_2\leq -z_3$, such that $\|x-\bar G(y)\|\geq \frac{1}{\|z_1\|} \frac{3z_3}{2}>\delta_1$. However, since $\|x-\bar G(y)\|_{P_s}^2<v_L$ implies $\|x-\bar G(y)\|<\delta_1$, a contradiction is obtained, proving that $\|x-\bar G(y)\|_{P_s}^2\leq v_L$ implies $\|\bar G(x)-y\|_{P_s}^2> v_L$.
\\
We have proven that, for all $x,y\in C$, the inequality $\|x-\bar G(y)\|_{P_s}^2\leq v_L$ implies the inequalities $\|x-y\|_{P_0}^2>v_L$ and $\|\bar G(x)-y\|_{P_s}^2>v_L$, such that $S_0\cap S_1\cap V_L=\emptyset$ and $S_1\cap S_2\cap V_L=\emptyset$. Interchanging the role of $x$ and $y$ in the reasoning above, we find $S_0\cap S_2\cap V_L=\emptyset$, which completes the proof of 1).
\\
To prove $S_1\cap V_L\cap (D\times (C\cup D))=\emptyset$, it suffices to observe that $\bpm{x^T&y^T}^T\in S_1\cap V_L$ implies $\|x-\bar G(y)\|<\delta_1$ and, by Lemma~\ref{lembarg}, $z_1 \bar G(y)+z_2>\tfrac{z_3}2$. Hence, $z_1x+z_2\geq z_1\bar G(y)+z_2-\|z_1\| \|x-\bar G(y)\|\geq \frac{z_3}{2}-\|z_1\| \delta_1>0$, which directly implies $x\not \in D$ due to \eqref{eqhsafd} and hence indeed $S_1\cap V_L\cap (D\times C\cup D)=\emptyset$. Since $V(x,y)=V(y,x)$, we also find that $S_2\cap V_L\cap ((C\cup D)\times D)=\emptyset$, such that 2) is proven.
\\
To prove 3), we observe that $\bpm{x^T&y^T}^T\in S_1\cap V_L$ yields $\|x-\bar G(y)\|_{P_s}^2\leq v_L$, which implies $\|x-\bar G(y)\|<\delta_1$, such that Lemma~\ref{lembarg} implies $z_1y+z_2\leq -z_3<0$. Similarly, we obtain $z_1x+z_2<0$ for $\bpm{x^T&y^T}^T\in S_2\cap V_L$, proving 3). This completes the proof of the lemma.
\end{pf}
The proof of Lemma~\ref{lemvdes} relies on Lemma~\ref{lemvcases} and the following result.
\begin{lem}\label{lemvbnds}
Consider the hybrid system \eqref{eqhsaf} and the Lyapunov function candidate $V$ in \eqref{eqv} satisfying the conditions of Lemma~\ref{lemvdes}. There exist $\mathcal{K}_\infty$-functions $\alpha_1,\alpha_2$ satisfying \eqref{eqvbnds}.
\end{lem}
\begin{pf}
To find lower and upper bounds on the function $V$, we introduce two scalars $\underline \lambda,\bar \lambda$, with
$\underline \lambda>0$ smaller than the eigenvalues of $P_0$ and $P_s$, and
$\bar \lambda>0$ larger than the eigenvalues of $P_0$ and $P_s$.  Let $e_0(x,y):=x-y$, $e_1(x,y):=x-\bar G(y)$ and $e_2(x,y):=\bar G(x)-y$, and define the nonnegative functions  $v_0(x,y):=\|e_0(x,y)\|_{P_0}^2$, $v_1(x,y):=\|e_1(x,y)\|_{P_s}^2$ and $v_2(x,y):=\|e_2(x,y)\|_{P_2}^2$.
\\
The definition of $\underline \lambda $ and $\bar \lambda$ yields
\beq{
\underline \lambda \|e_i(x,y)\|^2\leq v_i(x,y)\leq \bar \lambda \|e_i(x,y)\|^2,\ i=0,1,2.\label{ineqlams}
}
Observe that $V(x,y)=v_j(x,y)$, for some $j\in I$, with $I=\{0,1,2\}$.
The first inequality of \eqref{ineqlams} then yields $V(x,y)\geq \underline \lambda \|e_j\|^2 \geq \underline \lambda \min_{i\in I}\|e_i(x,y)\|^2$ holds. Note that $\min_{i\in I}\|e_i(x,y)\|^2\geq d(x,y)^2$, cf.\ \eqref{eqd}, where $\mathcal{A}=\{\bpm{z_x^T&z_y^T}^T\in C_{e}\cup D_{e}|\ x=y\lor x=G(y)\lor G(x)=y\}$ is used. Hence, we obtain $V(x,y)\geq \underline \lambda d(x,y)^2$. With $\alpha_1(s)=\underline \lambda s^2$, the first inequality in \eqref{eqvbnds} is satisfied.
\\
To derive an upper bound for the function $V$, we observe that $w\left(\bpm{x\\y}\right)=\sqrt{V(x,y)}$ is a Lipschitz function with Lipschitz constant $L_V=\sqrt{\bar \lambda}\sigma $, with $\sigma$ the maximum of the largest singular values of the matrices $\bpm{I_n&-I_n}$, $\bpm{I_n&-(L+MJ)}$ and $\bpm{I_n&-(L+MJ+sLJ^Tz_1)}$. Hence,
\[
w\left(\bpm{x\\y}\right)\leq w\left(\bpm{x^*\\y^*}\right)+L_V\left\|\bpm{x^*\\y^*}-\bpm{x\\y} \right\|
\]
for all $\bpm{x^*\\y^*}\in C_e\times D_e$. Given $x,y$, we can select $\bpm{x^*\\y^*}\in \arg \inf_{z\in \mathcal{A}}\left\|\bpm{x\\y}-z\right\|$, which, since $\mathcal{A}$ is non-empty and closed, cf.\ Theorem~\ref{thmdprop}, exists and is bounded. We observe that $V\left(\bpm{x^*\\y^*}\right)=0$ and, hence,  $w\left(\bpm{x^\star\\y^\star}\right)=0$. Consequently, we find
\bal{
w\left(\bpm{x\\y}\right)&\leq L_V\left\| \bpm{x\\y}- \arg \inf_{z\in \mathcal{A}}\left\|\bpm{x\\y}-z\right\|\nonumber
 \right\|,}
 \bal{
w\left(\bpm{x\\y}\right)&\leq L_V \inf_{z\in \mathcal{A}}\left\|\bpm{x\\y}-z\right\|= L_V d(x,y),}
where \eqref{eqd} is useed. Since $V(x,y)= w^2\left(\bpm{x\\y}\right)$, $\alpha_2(s)=L_V^2s^2$ satisfies \eqref{eqvbnds}, thereby proving the lemma.
\end{pf}
\begin{pf}[Proof of Lemma~\ref{lemvdes}]
To prove the lemma, first, we observe that Lemma~\ref{lemvbnds} directly guarantees that there exist functions $\alpha_1,\alpha_2$ satisfying \eqref{eqvbnds}. In addition, Lemma~\ref{lemvcases} directly proves that there exists a sufficiently small $v_L>0$ such that $V$ is smooth in an open domain containing $V_L$. It remains to be proven that \eqref{lyapj1}-\eqref{lyapj2} imply \eqref{eqvjump}.
\\
Recall the sets $S_0,S_1,S_2$ in \eqref{eqS012}. Jumps of \eqref{eqhsaf} may trigger jumps between these sets. From item 2) in Lemma~\ref{lemvcases}, we observe that for $(x,y)\in S_1\cap V_L$ and $(x,y)\in S_2\cap V_L$ jumps of $x$ and $y$, respectively, are not feasible. Consequently, when  $(x,y)\in S_0$, both $x$ and $y$ can jump, while from $(x,y)\in S_1$, only a jump of $y$ is feasible, and $(x,y)\in S_2$ implies $x\not \in D$. We will now prove that \eqref{eqvjump} holds along these four jumps:
\begin{itemize}
\item[a)]
We first study the jump $(x,y)\to (G(x),y)$, with $(x,y)\in S_0$. Since 3) of Lemma~\ref{lemvcases} implies that $\bar G(y)=(L+MJ)y+H+MK$ as $z_1y+z_2\leq 0$ and $x\in D$ implies $\bar G(x)=G(x)=(L+MJ)x+H+MK$, we observe that $V(G(x),y)\leq \|G(x)-\bar G(y)\|_{P_s}^2=\|\bar G(x)-\bar G(y)\|_{P_s}^2=(x-y)^T(L+MJ)^TP_s(L+MJ)(x-y)$, such that \eqref{lyapj1} implies that \eqref{eqvjump} holds.\footnote{
When $V(x,y)\leq e^{-\lambda_d}v_L$, we conclude that $\|G(x)-\bar G(y)\|_{P_s}^2<v_L$, which implies, by Lemma~\ref{lemvcases}, that $(G(x),y)\in S_1$, as shown in Fig.~\ref{fig3j}.}
\item[b)]
For a jump $(x,y)\to (x,G(y))$ with $(x,y)\in S_1$, we observe $V(x,G(y))\leq \|x-G(y)\|^2_{P_0}=\|x-\bar G(y)\|_{P_0}^2$, as $y\in D$. Hence, \eqref{lyapj2} implies \eqref{eqvjump} in this case.\footnote{
Again, when before the jump, $V(x,y)\leq e^{-\lambda_d}v_L$ holds, we conclude that $\|x-G(y)\|_{P_0}^2\leq v_L$, such that $(x,G(y))\in S_0$, as shown in Fig.~\ref{fig3j}, follows.}
\item[c)] For a jump $(x,y)\to (x,G(y))$, with $(x,y)\in S_0$, \eqref{eqvjump} directly follows from combining a) with the symmetry relation $V(x,y)=V(y,x)$.
\item[d)] For a jump $(x,y)\to (G(x),y)$ with $(x,y)\in S_2$, symmetry of $V$ and b) imply \eqref{eqvjump}.
\end{itemize}
Hence, we have proven that \eqref{eqvjump} holds over all feasible jumps, therewith concluding the proof of the lemma.
\end{pf}
\begin{pf}[Proof of Theorem~\ref{thmtrack}]
We prove this theorem by application of Theorem \ref{thmgen}. Lemma~\ref{lemvdes} proves that \eqref{eqvbnds} and \eqref{eqvjump} hold for some $v_L>0$. Hence, we now show that the assumptions in the theorem prove that \eqref{eqvflow} is satisfied in the sub-level set $V_L=V^{-1}([0,v_L])$.
\\
According to Lemma~\ref{lemvdes}, $V$ is differentiable in $V_L$, such that we evaluate $ \langle \left.\frac{\partial V}{\partial q}\right|_q,f \rangle$ for $f\in F_e(t,\bar x_d(t,j),y)$ only when $q=\bpm{\bar x_d(t,j)^T&y^T}^T\in V_L\cap C_{e}$, where, for almost all $t$, $F_e$ is single-valued, and we distinguish the three cases given by the minimisers of \eqref{eqv}.
If $(\bar x_d(t,j),y)\in S_0\cap V_L$, then
\[
\frac{\partial V}{\partial q}=2(\bar x_d(t,j)-y)^TP_0\bpm{I_n&-I_n}\]
and
\[
F_e=\bpm{A\bar x_d(t,j)+E+Bu_\text{ff}(t)\\Ay+E+B(u_\text{ff}(t)-c_0(\bar x_d(t,j)-y))},
\]
such that \eqref{eqvflow} is guaranteed by \eqref{lmi1}.
\\
If $(\bar x_d(t,j),y)\in S_1\cap V_L$, then 3) of Lemma~\ref{lemvcases} implies $\bar G(y)=(L+MJ)y+H+MK$. Consequently
\[
\tfrac{\partial V}{\partial q}=2s_1^TP_s\bpm{I_n&-(L+MJ)}
\]
and
\[
F_e(t,\bar x_d(t,j),y)\!=\!\bpm{A\bar x_d(t,j)+E+Bu_\text{ff}(t)\\Ay\!+\!E\!+\!B(u_\text{ff}(t)
\!-\!\tfrac{\beta_2^T\beta_1(t)}{\beta_2^T\beta_2}\!-\!c_1s_1)}
\]
 with $s_1=\bar x_d(t,j)-\bar G (y)$ holds. Hence, we obtain
$\frac{\partial V}{\partial q}F_e(t,x,y)=2s_1^TP_s(A\bar x_d(t,j)+E+Bu_{\text{ff}}(t)-(L+MJ)Ay-(L+MJ)E-(L+MJ)Bu_{\text{ff}(t)}-\tfrac{\beta_2\beta_2^T}{\beta_2^T\beta_2}\beta_1(t)-(L+MJ)Bc_1s_1.$
With \eqref{eqbetspans}, we find $\tfrac{\beta_2\beta_2^T}{\beta_2^T\beta_2}\beta_1(t)=\beta_1(t)$, such that
\bal{\tfrac{\partial V}{\partial q}&F_e(t,\bar x_d(t,j),y)=\nonumber \\
&
2s_1^TP_s(A\bar x_d(t,j)+(I-L-MJ)(E+Bu_{\text{ff}}(t))\nonumber \\&-(L+MJ)Ay-\beta_1(t)+\beta_2c_1s_1).
}
Since $y=(L+MJ)^{-1}(-s_1+\bar x_d(t,j)-H-MK)=-(L+MJ)^{-1}s_1+\bar G^{\circ}(\bar x_d(t,j))$, (which holds as $\max(0,z_1\bar x_d(t,j)+z_2)=0$ follows from $(\bar x_d(t,j),y)\in S_1\cap V_L$, cf.\ Lemma~\ref{lemvcases}) we obtain
\bal{\tfrac{\partial V}{\partial q}&F_e(t,\bar x_d(t,j),y)=\nonumber \\&
2s_1^TP_s((L+MJ)A(L+MJ)^{-1}+\beta_2c_1)s_1,
}
where we used the design of $\beta_1$. Hence, \eqref{lmi3} guarantees that \eqref{eqvflow} holds in this case.
\\
Now, we focus on the case $(\bar x_d(t,j),y)\in S_2\cap V_L$. In that case, from  3) of Lemma~\ref{lemvcases}, we observe that $\max(0,z_1y+z_2)=0$ follows from $(\bar x_d(t,j),y)\in S_2\cap V_L$, cf.\ Lemma~\ref{lemvcases}. Hence,
\[
\frac{\partial V}{\partial q}=2s_2^TP_s\bpm{(L+MJ)&-I_n}
\] and
\[
F_e(t,\bar x_d(t,j),y)=\bpm{A\bar x_d(t,j)\!+\!E\!+\!Bu_\text{ff}(t)\\Ay\!+\!E\!+\!B(u_\text{ff}(t)
\!-\!\frac{\beta_4^T\beta_3(t)}{\beta_4^T\beta_4}\!+\!c_2s_2)}
\] with $s_2=\bar G(\bar x_d(t,j))\!-\!y$.
From \eqref{eqbetspans} follows
$\frac{\beta_4\beta_4^T}{\beta_4^T\beta_4}\beta_3(t)=\beta_3(t)$, such that
$\frac{\partial V}{\partial q}F_e(t,\bar x_d(t,j),y)=2s_2^TP_s\big(As_2-\beta_3(t)+\bpm{L\!\!+\!\!MJ&-I_{n}}$ $\bpm{A\bar x_d(t,j)+Bu_{\text{ff}}(t)+E\\A\bar G(\bar x_d(t,j))+Bu_{\text{ff}}(t)+E}+Bc_2s_2\big)$, where we used
$y=\bar G(\bar x_d(t,j))-s_2$. With the design of $\beta_3, \beta_4$, we find
\bal{
\tfrac{\partial V}{\partial q}F_e(t,\bar x_d(t,j),y)=2s_2^TP_s(A+\beta_4c_2)s_2,
}
such that \eqref{lmi2} proves that \eqref{eqvflow} holds in this case.
Consequently, if \eqref{lmi1}-\eqref{lmi2} hold, \eqref{eqvflow} is obtained. Hence, Theorem~\ref{thmgen} proves that $x_d$ is asymptotically stable with respect to $d$.
\end{pf}

\end{document}